\documentclass[11pt]{article}
\usepackage{latexsym,amssymb}
\usepackage{amsmath}
\usepackage[mathscr]{eucal}
\usepackage{color}
\usepackage[abbrev]{amsrefs}

\newtheorem{Theorem}{\bf Theorem}[section]
\newtheorem{Lemma}{\bf Lemma}[section]
\newtheorem{Proposition}{\bf Proposition}[section]  
\newtheorem{Corollary}{\bf Corollary}[section]
\newtheorem{Remark}{\bf Remark}[section]
\newtheorem{Example}{\bf Example}[section]
\newtheorem{Definition}{\bf Definition}[section]

\newenvironment{theorem}{\begin{Theorem}$\!\!\!$}{\end{Theorem}}
\newenvironment{lemma}{\begin{Lemma}$\!\!\!$}{\end{Lemma}}
\newenvironment{proposition}{\begin{Proposition}$\!\!\!$}{\end{Proposition}}
\newenvironment{corollary}{\begin{Corollary}$\!\!\!$}{\end{Corollary}}
\newenvironment{remark}{\begin{Remark}$\!\!\!$}{\end{Remark}}

\newenvironment{definition}{\begin{Definition}$\!\!\!$}{\end{Definition}}

\newcommand{\R}{{\bf{R}}}

\numberwithin{equation}{section}

\def\XXint#1#2#3{{\setbox0=\hbox{$#1{#2#3}{\int}$}
\vcenter{\hbox{$#2#3$}}\kern-.5\wd0}}

\begin{document}
\title{Decay estimates for Schr\"odinger heat semigroup\\ 
with inverse square potential in Lorentz spaces}
\author{Kazuhiro Ishige and Yujiro Tateishi\\
Graduate School of Mathematical Sciences, The University of Tokyo\\
3-8-1 Komaba, Meguro-ku, Tokyo 153-8914, Japan\vspace{10pt}
}
\date{}
\maketitle
\noindent
\qquad\qquad
E-mail addresses:

\noindent
\qquad\qquad
ishige@ms.u-tokyo.ac.jp (K. Ishige), 
tateishi@ms.u-tokyo.ac.jp (Y. Tateishi)
$$
\qquad
$$
\begin{abstract}
Let $H:=-\Delta+V$ be a nonnegative Schr\"odinger operator on $L^2({\bf R}^N)$, 
where $N\ge 2$ and $V$ is an inverse square potential. 
In this paper we obtain sharp decay estimates of the operator norms of $e^{-tH}$ and $\nabla e^{-tH}$ 
in Lorentz spaces.
\end{abstract}

%Notation:
%\begin{itemize}
%  \item Page 3: $A_\lambda^\pm$, $D_\lambda$, $A_{1,k}$, $A_{2.k}$, $B_k$, $v_k$
%  \item 
%  \item 
%\end{itemize}
\newpage
%%%%%%%%%%%%%%%%%%%%%%%%%%%%%%%%%%%%
%%%%%%%%%%%%%%%%%%%%%%%%%%%%%%%%%%%%
\section{Introduction}
%%%%%%%%%%%%%%%%%%%%%%%%%%%%%%%%%%%%
%%%%%%%%%%%%%%%%%%%%%%%%%%%%%%%%%%%%
Nonnegative Schr\"odinger operators $H:=-\Delta+V$ on $L^2({\bf R}^N)$ and their heat semigroups $e^{-tH}$ 
have been studied by many mathematicians since the pioneering work due to Simon~\cite{S}. 
See e.g. \cite{Bar}, \cite{CK}, \cite{CUR}, \cite{DS}, 
\cite{IIY01}, \cite{IIY02}, \cite{IK01}--\cite{IKO}, 
\cite{LS}--\cite{M}, 
\cite{P3}--\cite{Zhang} and references therein.   
(See also the monographs of Davies \cite{Dav}, Grigor'yan \cite{Gri} and Ouhabaz \cite{Ouh}.) 
Decay estimates of operator norms of $\nabla^\ell e^{-tH}$ 
are fundamental to the study of nonnegative Schr\"odinger operators and their related fields.
Here $\nabla:=(\partial/\partial x_1,\dots,\partial/\partial x_N)$ and $\ell\in\{0,1,\dots\}$.
However, the large time decay estimates of operator norms of $\nabla^\ell e^{-tH}$ 
are delicate and they are widely open even in Lebesgue spaces 
(see e.g. \cite{DS}, \cite{Ishige}--\cite{IK05}, \cite{LL}, \cite{S} and references therein). 

In this paper we focus on a nonnegative Schr\"odinger operator $H:=-\Delta+V$, 
where $V$ is an inverse square potential, more precisely, $V$ satisfies 
the following condition~($\mbox{V}$): 
\begin{equation}
\tag{$\mbox{V}$}
\left\{
\begin{array}{ll}
({\rm i}) & \mbox{$V=V(|x|)$ in ${\bf R}^N\setminus\{0\}$ and $V\in C^1((0,\infty))$};\vspace{7pt}\\
({\rm ii}) & \mbox{$V(r)=\lambda_1r^{-2}+O(r^{-2+\rho_1})$ as $r\to +0$},\vspace{3pt}\\
 & \mbox{$V(r)=\lambda_2r^{-2}+O(r^{-2-\rho_2})$ as $r\to \infty$},\vspace{3pt}\\
  & \mbox{for some $\lambda_1$, $\lambda_2\in[\lambda_*,\infty)$ with $\lambda_*:=-(N-2)^2/4$ and $\rho_1$, $\rho_2>0$};\vspace{7pt}\\
({\rm iii}) &  \displaystyle{\sup_{r>0}\,\left|\,r^3 \frac{d}{dr}V(r)\right|<\infty}.
\end{array}
\right.
\end{equation}
The purpose of this paper is to obtain the decay estimates of 
$$
\|\nabla^\ell e^{-tH}\|_{(L^{p,\sigma}\to L^{q,\theta})}
:=\sup
\left\{\left\|\nabla^\ell e^{-tH}\phi\right\|_{L^{q,\theta}({\bf R}^N)}\,:\,
\mbox{$\phi\in C_{\rm c}({\bf R}^N)$ with $\|\phi\|_{L^{p,\sigma}({\bf R}^N)}=1$}\right\},
$$
where $\ell\in\{0,1\}$ and 
$$
(p,q,\sigma,\theta)\in\Lambda:=
\left\{1\le p\le q\le\infty,\,\sigma,\,\theta\in[1,\infty]:
\begin{array}{ll}
\mbox{$\sigma=1$\, if \,$p=1$}, & \mbox{$\sigma=\infty$\, if \,$p=\infty$}\vspace{3pt}\\
\mbox{$\theta=1$\, if \,$q=1$}, & \mbox{$\theta=\infty$\, if \,$q=\infty$}\vspace{3pt}\\
\mbox{$\sigma\le\theta$\, if \,$p=q$} &
\end{array}
\right\}.
$$
Here $e^{-tH}\phi$ is a solution to the Cauchy problem
\begin{equation}
\tag{P}
\left\{
\begin{array}{ll}
\partial_t u=\Delta u-V(|x|)u & \quad\mbox{in}\quad{\bf R}^N\times(0,\infty),\vspace{5pt}\\
u(\cdot,0)=\phi  & \quad\mbox{in}\quad{\bf R}^N,
\end{array}
\right.
\end{equation}
and $\|\nabla^\ell e^{-tH}\|_{(L^{p,\sigma}\to L^{q,\theta})}$ is 
the operator norm of $\nabla^\ell e^{-tH}$ 
from the Lorentz space $L^{p,\sigma}({\bf R}^N)$ to $L^{q,\theta}({\bf R}^N)$. 
This paper can be regarded as a continuation of \cite{IIY02}, 
where the large time sharp decay estimates of $\|e^{-tH}\|_{(L^{p,\sigma}\to L^{q,\theta})}$ were obtained when $V\in C^1([0,\infty))$ 
and $\lambda_2>\lambda_*$. 
\vspace{3pt}

We introduce some notations. 
For $x\in{\bf R}^N$ and $R>0$, 
let $B(x,R):=\{y\in{\bf R}^N\,:\,|y-x|<R\}$ and $B(x,R)^c:={\bf R}^N\setminus B(x,R)$.  
For any $r\in[1,\infty]$, let $r'$ be the H\"older conjugate number of $r$, that is, 
$$
r'=\frac{r}{r-1}\quad\mbox{if}\quad 1<r<\infty,
\quad
r'=1\quad\mbox{if}\quad r=\infty,
\quad
r'=\infty\quad\mbox{if}\quad r=1.
$$
Let $\Delta_{{\bf S}^{N-1}}$ be the Laplace-Beltrami operator on ${\bf S}^{N-1}$. 
Let $\{\omega_k\}_{k=0}^\infty$ be 
the eigenvalues of 
\begin{equation}
\label{eq:1.1}
-\Delta_{{\bf S}^{N-1}}Q=\omega Q\quad\mbox{on}\quad{\bf S}^{N-1},
\qquad
Q\in L^2({\bf S}^{N-1}). 
\end{equation}
Then $\omega_k=k(N+k-2)$ for $k\in\{0,1,2,\dots\}$. 
Let
$\{Q_{k,i}\}_{i=1}^{d_k}$ and $d_k$ be 
the orthonormal system and the dimension of the eigenspace corresponding to $\omega_k$, respectively. 
Here 
\begin{equation}
\label{eq:1.2}
d_k=\frac{(N+2k-2)(N+k-3)!}{(N-2)!k!}=O(k^{N-2})\quad\mbox{as}\quad k\to\infty.
\end{equation}

Assume condition~($\mbox{V}$) 
and let $H:=-\Delta+V$ be nonnegative, that is, 
\begin{equation}
\tag{N}
\int_{{\bf R}^N}\left[|\nabla\phi|^2+V(|x|)\phi^2\right]\,dx\ge 0,
\qquad\phi\in C_c^\infty({\bf R}^N\setminus\{0\}). 
\end{equation}
The operator $H$ is said {\it subcritical} if, for any $W\in C_c({\bf R}^N)$, 
$H-\epsilon W$ is nonnegative for small enough $\epsilon>0$. 
If not, the operator $H$ is said {\it critical}. 

For any $k\in\{0,1,2,\dots\}$, set $A_{1,k}:=A^+_{\lambda_1+\omega_k}$ and 
\begin{equation*}
\begin{split}
A_{2,k}:= & 
\left\{
\begin{array}{ll}
A^-_{\lambda_2} & \mbox{if $k=0$ and $H$ is critical},\vspace{3pt}\\
A^+_{\lambda_2+\omega_k} & \mbox{otherwise},
\end{array}
\right.\\
B_k:= & 
\left\{
\begin{array}{ll}
1 &  \mbox{if $k=0$, $\lambda_2=\lambda_*$ and $H$ is subcritical},\vspace{3pt}\\
0 &  \mbox{otherwise}.
\end{array}
\right.
\end{split}
\end{equation*}
Here 
\begin{equation}
\label{eq:1.3}
A^\pm_\lambda:=\frac{-(N-2)\pm\sqrt{D_\lambda}}{2}
\quad\mbox{for $\lambda\ge\lambda_*$, 
where $D_\lambda:=(N-2)^2+4\lambda$}.
\end{equation}
By the standard theory for ordinary differential equations 
we see that, 
for any $k\in\{0,1,2,\dots\}$, 
there exists a unique solution $h_k=h_k(r)$ to the problem
\begin{equation}
\label{eq:1.4}
\begin{split}
 & \frac{d^2}{dr^2}h_k+\frac{N-1}{r}\frac{d}{dr}h_k-V_k(r)h_k=0\quad\mbox{in}\quad(0,\infty),\\
 & h_k(r)=r^{A_{1,k}}(1+o(1))\quad\mbox{as}\quad r\to+0,
\end{split}
\end{equation}
where $V_k(r):=V(r)+\omega_kr^{-2}$. 
(See also Section~2.1.) 
Notice that $h_k\in L^2(B(0,1))$. 
Furthermore, it follows from \cite[Theorem~1.1]{IKO} that
\begin{equation}
\label{eq:1.5}
h_k(r)=c_kv_k(r)(1+o(1))\quad\mbox{as}\quad r\to\infty,
\quad\mbox{where}\quad
v_k(r):=r^{A_{2,k}}(\log r)^{B_k},
\end{equation}
for some $c_k>0$.
The function~$h_0$ is said {\it a positive harmonic function} for the operator $H$
and it plays an important role in our analysis. 
When $H$ is critical, if $h_0\not\in L^2({\bf R}^N)$, 
then $H$ is said {\it null-critical}\,: if not, $H$ is said {\it positive-critical}. 
The decay of the fundamental solution $p=p(x,y,t)$ corresponding to $e^{-tH}$ 
depends on whether $H$ is either subcritical, null-critical or positive critical. 
See \cite{P3}. (See also \cite{IKO}.)
%\begin{itemize}
%  \item 
%  Let $H$ be subcritical and $x$, $y\in{\bf R}^N\setminus\{0\}$ with $x\not=y$. Then
%  $$
%  \lim_{t\to\infty} p(x,y,t)=0,\qquad\int_0^\infty p(x,y,t)\,dt<\infty.
%  $$
%  \item 
%  Let $H$ be critical and $x$, $y\in{\bf R}^N\setminus\{0\}$ with $x\not=y$. 
%  If $H$ is null-critical, that is, $A_{2,0}\ge -N/2$, then 
%  $$
%  \lim_{t\to\infty} p(x,y,t)=0,\qquad\int_0^\infty p(x,y,t)\,dt=\infty.
%  $$
%  If $H$ is positive-critical, that is, $A_{2,0}<-N/2$, then 
%  $$
%  \lim_{t\to\infty} p(x,y,t)=\frac{h_0(|x|)h_0(|y|)}{\|h_0\|_{L^2({\bf R}^N)}^2}.
%  $$
%\end{itemize}
%%
In this paper 
we assume either
\begin{equation}
\tag{N'}
\qquad\qquad
{\rm (i)}\quad\mbox{$H$ is subcritical}\quad\mbox{or}\quad
{\rm (ii)}\quad\mbox{$H$ is critical and $A_{2,0}>-N/2$},
\qquad
\end{equation}
and obtain decay estimates of $\|\nabla^\ell e^{-tH}\|_{(L^{p,\sigma}\to L^{q,\theta})}$ as $t\to +0$ and as $t\to\infty$, 
where $\ell\in\{0,1\}$ and $(p,q,\sigma,\theta)\in\Lambda$. 
Case~(ii) is in the null-critical one. 
We remark that $e^{-tH}\phi$ does not necessarily decay as $t\to\infty$ if $H$ is positive-critical. 
(See also Remark~\ref{Remark:1.1}~(iii).)
\vspace{5pt}

Now we are ready to state the main results of this paper. 
\begin{theorem}
\label{Theorem:1.1}
Assume conditions~{\rm ($\mbox{V}$)} and {\rm (N')}. 
Let $(p,q,\sigma,\theta)\in\Lambda$, $j\in\{0,1,2,\dots\}$ and $\ell\in\{0,1\}$. 
Then there exists $C>0$ such that
\begin{equation}
\label{eq:1.7}
\left\|\partial_t^j\nabla^\ell e^{-tH}\right\|_{(L^{p,\sigma}\to L^{q,\theta})}
\le Ct^{-\frac{N}{2}-j}\frac{\|h_0\|_{L^{p',\sigma'}(B(0,\sqrt{t}))}}{h_0(\sqrt{t})}
\left(\frac{\|\nabla^\ell h_0\|_{L^{q,\theta}(B(0,\sqrt{t}))}}{h_0(\sqrt{t})}+t^{\frac{N}{2q}-\frac{\ell}{2}}\right)
\end{equation}
for $t>0$.
\end{theorem}
We remark that the right-hand side of inequality~\eqref{eq:1.7} possibly diverges for some $(p,q,\sigma,\theta)\in\Lambda$. 
As a corollary of Theorem~\ref{Theorem:1.1}, 
we have:
\begin{corollary}
\label{Corollary:1.1}
Assume conditions~{\rm ($\mbox{V}$)} and {\rm (N')}. 
Let $(p,q,\sigma,\theta)\in\Lambda$. 
Then there exists $C>0$ such that
\begin{equation}
\label{eq:1.7a}
\left\|e^{-tH}\right\|_{(L^{p,\sigma}\to L^{q,\theta})}
\le Ct^{-\frac{N}{2}}\frac{\|h_0\|_{L^{p',\sigma'}(B(0,\sqrt{t}))}\|h_0\|_{L^{q,\theta}(B(0,\sqrt{t}))}}{h_0(\sqrt{t})^2}
\quad\mbox{for}\quad t>0.
\end{equation}
\end{corollary}
\begin{remark}
\label{Remark:1.1}
{\rm (i)} 
Assume conditions~{\rm ($\mbox{V}$)} and {\rm (N')}. 
Furthermore, assume that $V\in C^1([0,\infty))$ and $\lambda_2>\lambda_*$. 
Then the large time sharp decay estimates of 
$\left\|e^{-tH}\right\|_{(L^{p,\sigma}\to L^{q,\theta})}$ have been already obtained in \cite[Theorem~1.1]{IIY02}. 
Our decay estimate \eqref{eq:1.7a} gives the same decay estimate as in \cite[Theorem~1.1]{IIY02} 
and it has a simpler expression. 
\vspace{3pt}
\newline
{\rm (ii)} 
In a forthcoming paper~\cite{IT02},
under conditions {\rm (V)} and {\rm (N')}, 
we obtain decay estimates of $\left\|\nabla^\ell e^{-tH}\right\|_{(L^{p,\sigma}\to L^{q,\theta})}$ for $\ell\in\{0,1,2,\dots\}$. 
Furthermore, we find $C>0$ such that
$$
\left\|\nabla^\ell e^{-tH}\right\|_{(L^{p,\sigma}\to L^{q,\theta})}
\ge C^{-1}t^{-\frac{N}{2}}\frac{\|h_0\|_{L^{p',\sigma'}(B(0,\sqrt{t}))}}{h_0(\sqrt{t})}
\left(\frac{\|\nabla^\ell h_0\|_{L^{q,\theta}(B(0,\sqrt{t}))}}{h_0(\sqrt{t})}+t^{\frac{N}{2q}-\frac{\ell}{2}}\right)
$$
for $t>0$, where $\ell\in\{0,1\}$, 
and show the sharpness of decay estimate~\eqref{eq:1.7}.
\vspace{3pt}
\newline
{\rm (iii)} Under condition~{\rm (V)}, 
case~{\rm (ii)} covers all of the null-critical cases except for the case when $A_{2,0}=-N/2$. 
On the other hand, 
in this paper, we don't treat the case when $A_{2,0}=-N/2$ 
since it is on the borderline between null-criticality and positive-criticality
and it is too delicate. 
\end{remark}

We prove Theorem~\ref{Theorem:1.1} by developing the arguments in \cite{IK01} 
and combining the results of \cite{IKO} with parabolic regularity theorems. 
For any $\phi\in C_c({\bf R}^N)$, 
we find radially symmetric functions 
$\{\phi_{k,i}\}_{k=0,1,\dots,\,i=1,\dots,d_k}\subset L^2({\bf R}^N)$ 
such that 
\begin{equation}
\label{eq:1.8}
\phi(x)=\sum_{k=0}^\infty\sum_{i=1}^{d_k}\phi_{k,i}(|x|)Q_{k,i}\left(\frac{x}{|x|}\right)
\quad\mbox{in}\quad L^2({\bf R}^N).
\end{equation}
Let $H_k:=-\Delta+V_k(|x|)$ and set 
\begin{equation}
\label{eq:1.9}
v_{k,i}(|x|,t):=[e^{-tH_k}\phi_{k,i}](|x|),
\qquad
u_{k,i}(x,t)=v_{k,i}(|x|,t)Q_{k,i}\left(\frac{x}{|x|}\right).
\end{equation}
Then 
\begin{equation}
\label{eq:1.10}
\left[e^{-tH}\phi\right](x)=\sum_{k=0}^\infty\sum_{i=1}^{d_k} u_{k,i}(x,t)
=\sum_{k=0}^\infty\sum_{i=1}^{d_k} v_{k,i}(|x|,t)Q_{k,i}\left(\frac{x}{|x|}\right)
\quad\mbox{in}\quad C^2(K)
\end{equation}
for compact sets $K\subset{\bf R}^N\setminus\{0\}$ and $t>0$ 
(see \cite{IM} and \cite{IKM}). 
For the proof of Theorem~\ref{Theorem:1.1}, 
we obtain uniform estimates of $\{v_{k,i}\}$ inside parabolic cones 
by constructing supersolutions. 
This is the main difficulty for the proof of Theorem~\ref{Theorem:1.1} and requires delicate analysis. 
Then, due to the radial symmetry of $v_{k,i}$ and parabolic regularity theorems, 
we obtain decay estimates of $\nabla^\ell e^{-tH}\phi$ inside parabolic cones. 
On the other hand, we obtain estimates of $\nabla^\ell e^{-tH}\phi$ outside parabolic cones 
by the Gaussian estimate of the fundamental solution $p=p(x,y,t)$ and parabolic regularity theorems. 
Combining the estimates of $\nabla^\ell e^{-tH}\phi$ inside and outside parabolic cones, 
we obtain Theorem~\ref{Theorem:1.1}.

The rest of this paper is organized as follows. 
In Section~2 we formulate a definition of the solution to problem~(P) 
and Lorentz spaces. 
Furthermore, we study the asymptotic behavior of $h_k$ as $r\to +0$ and as $r\to\infty$. 
In Section~3 we obtain upper estimates of $e^{-tH}\phi$ by using the Gaussian estimate of $p=p(x,y,t)$. 
In Section~4, combining the estimates in Sections~2 and~3, 
we obtain uniform estimates of $\{v_{k,i}\}$ inside parabolic cones with respect to $k$ and $i$.
In Section~5 we complete the proof of Theorem~\ref{Theorem:1.1}. 
Furthermore, we prove Corollary~\ref{Corollary:1.1}.
%%%%%%%%%%%%%%%%%%%%%%%%%%%%%%%%%%%%
%%%%%%%%%%%%%%%%%%%%%%%%%%%%%%%%%%%%
\section{Preliminaries}
%%%%%%%%%%%%%%%%%%%%%%%%%%%%%%%%%%%%
%%%%%%%%%%%%%%%%%%%%%%%%%%%%%%%%%%%%
In this section we formulate a definition of the solution to problem~(P) and introduce Lorentz spaces. 
We also obtain uniform estimates of the solution~$h_k$ to problem~\eqref{eq:1.4} with respect to~$k$. 

Throughout this paper, for any positive functions $f$ and $g$ on a set $E$, 
we write $f\asymp g$ for $x\in E$ if there exists $c>0$ such that 
$c^{-1}\le f(x)/g(x)\le c$ for $x\in E$.
By the letter $C$ we denote generic positive constants 
and they may have different values also within the same line. 
%%%%%%%%%%%%
\subsection{Definition of solutions to problem~($\mbox{P}$) and Lorentz spaces}
%%%%%%%%%%%%
Before considering problem~($\mbox{P}$), 
we formulate a definition of the solution to the problem
\begin{equation}
\tag{$\mbox{W}$}
\left\{
\begin{array}{ll}
\partial_t w=\displaystyle{\frac{1}{\nu}\mbox{div}\,(\nu\nabla w)} & 
\quad\mbox{in}\quad{\bf R}^N\times(0,\infty),\vspace{5pt}\\
w(\cdot,0)=\phi_*  & \quad\mbox{in}\quad{\bf R}^N,
\end{array}
\right.
\end{equation}
where $\nu:=h_0^2$ and $\phi_*\in L^2({\bf R}^N,\nu\,dx)$. 
\begin{definition}
\label{Definition:2.1}
Let $\phi_*\in L^2({\bf R}^N,\nu\,dx)$. 
A measurable function $w$ in ${\bf R}^N\times(0,\infty)$ 
is said a solution to problem~{\rm ($\mbox{W}$)} if 
$$
w\in L^\infty(0,\infty:L^2({\bf R}^N,\nu\,dx))\cap L^2(0,\infty:H^1({\bf R}^N,\nu\,dx))
$$
and $w$ satisfies
$$
\int_{{\bf R}^N}\psi(x,0)\phi_*(x)\nu(x)\,dx+
\int_0^\infty\int_{{\bf R}^N}
\left\{-w\partial_t\psi+\nabla w\nabla\psi\right\}\nu(x)\,dx\,dt=0
$$
for $\psi\in C^\infty_{{\rm c}}({\bf R}^N\times[0,\infty))$.
\end{definition}
Problem~($\mbox{W}$) possesses a unique solution $w$ such that 
$$
\|w(t)\|_{L^2({\bf R}^N,\nu\,dx)}\le\|\phi_*\|_{L^2({\bf R}^N,\nu\,dx)}\quad\mbox{for}\quad t>0. 
$$
See \cite[Section~2.1]{IM}. 
We denote by $e^{-tH_*}\phi_*$ the unique solution to problem~($\mbox{W}$). 
%We define the solution to problem~($\mbox{P}$) by the use of problem~($\mbox{W}$). 
%
\begin{definition}
\label{Definition:2.2}
Let $\phi\in L^2({\bf R}^N)$. 
Set 
$$
\left[e^{-tH}\phi\right](x):=h_0(|x|)\left[e^{-tH_*}\phi_*\right](x)
\quad\mbox{with}\quad
\phi_*(x):=\frac{\phi(x)}{h_0(|x|)}
$$
for $x\in{\bf R}^N$ and $t>0$. 
Then the function~$e^{-tH}\phi$ is said a solution to problem~{\rm ($\mbox{P}$)}.
\end{definition}
Notice that $\phi\in L^2({\bf R}^N)$ if and only if $\phi_*\in L^2({\bf R}^N,\nu\,dx)$. 
\vspace{3pt}

We define the Lorentz spaces. 
For any measurable function $\phi$ in ${\bf R}^N$, 
we denote by $\mu=\mu(\lambda)$ the distribution function of $\phi$, that is, 
$$
\mu(\lambda):=\left|\{x\in{\bf R}^N\,:\,|\phi(x)|>\lambda\}\right| 
\quad\mbox{for}\quad \lambda > 0.
$$
Here $|E|$ is the $N$-dimensional Lebesgue measure of $E$ for measurable sets $E$ in ${\bf R}^N$.
We define 
the non-increasing rearrangement $\phi^*$ of $\phi$
and the spherical rearrangement $\phi^{\sharp}$ of $\phi$ by 
$$
\phi^{*}(s):=\inf\{\lambda>0\,:\,\mu(\lambda)\le s\},
\qquad
\phi^\sharp(x):=\phi^*(\alpha_N|x|^N),
$$
for $s>0$ and $x\in {\bf R}^N$, respectively,
where $\alpha_N$ is the volume of the unit ball in ${\bf R}^N$. 
For any $(p,p,\sigma,\sigma)\in\Lambda$, 
we define the Lorentz space $L^{p,\sigma}({\bf R}^N)$ by 
$$
L^{p,\sigma}({\bf R}^N):=\{\phi\,:\, \mbox{$\phi$ is measurable in ${\bf R}^N$},\,\,\, \|\phi\|_{L^{p,\sigma}}<\infty\},
$$
where 
$$
\|\phi\|_{L^{p,\sigma}}:=
 \left\{
\begin{array}{ll}
  \displaystyle{\biggr(
   \int_{{\bf R}^N}\left(|x|^{\frac{N}{p}}\phi^{\sharp}(x)\right)^{\sigma}\frac{dx}{|x|^N}
  \biggr)^{\frac{1}{\sigma}}}\quad & \mbox{if}\quad 1\le \sigma<\infty, \vspace{5pt}\\
 \displaystyle{\sup_{x\in {\bf R}^N}\,|x|^{\frac{N}{p}}\phi^{\sharp}(x)}\quad & \mbox{if}\quad\sigma=\infty.
\end{array}
 \right.
$$
Here $N/p=0$ if $p=\infty$. 
The Lorentz spaces have the following properties:  
\begin{equation*}
\begin{array}{ll}
L^{p,p}({\bf R}^N)=L^p({\bf R}^N) & \mbox{if $1\le p\le \infty$};\vspace{3pt}\\ 
L^{p,\sigma_1}({\bf R}^N)\subset L^{p,\sigma_2}({\bf R}^N) & \mbox{if $1\le p<\infty$ and $1\le \sigma_1\le \sigma_2 \le \infty$}.\vspace{3pt}\\ 
\end{array}
\end{equation*}
Furthermore, 
there exists $C>0$ depending only on $N$ such that 
\begin{eqnarray}
\nonumber
 & & \|f+g\|_{L^{p,\sigma}}\le C(\|f\|_{L^{p.\sigma}}+\|g\|_{L^{p,\sigma}})\quad\mbox{if $f,g\in L^{p,\sigma}({\bf R}^N)$},\\
\label{eq:2.1}
 & & \|fg\|_{L^1}\le C\|f\|_{L^{p,\sigma}}\|g\|_{L^{p',\sigma'}}
 \qquad\qquad\,\,\mbox{if $f\in L^{p,\sigma}({\bf R}^N)$, $g\in L^{p',\sigma'}({\bf R}^N)$},\\
\label{eq:2.2}
 & & \|f*g\|_{L^{q,\theta}}\le C\|f\|_{L^{p,\sigma}}\|g\|_{L^{r,s}}
 \qquad\quad\,\,\,\mbox{if $f\in L^{p,\sigma}({\bf R}^N)$, $g\in L^{r,s}({\bf R}^N)$}.
\end{eqnarray}
Here $(q,q,\theta,\theta)$, $(r,r,s,s)\in\Lambda$ and
\begin{equation*}
   \frac{1}{r}+\frac{1}{p}=\frac{1}{q}+1,\qquad \frac{1}{\theta}=\frac{1}{s}+\frac{1}{\sigma}.
\end{equation*}
(See e.g. \cite{BS} and \cite{Gra}.) 
For any measurable function~$f$ in a domain $\Omega$, 
we say that $f\in L^{p,\sigma}(\Omega)$ if and only if $\tilde{f}\in L^{p,\sigma}({\bf R}^N)$, 
where $\tilde{f}$ is the zero extension of $f$ to ${\bf R}^N$. 
Furthermore, we write $\|f\|_{L^{p,\sigma}(\Omega)}=\|\tilde{f}\|_{L^{p,\sigma}}$.
Then, for any $R>0$, 
the function $f_A$ defined by $f_A(x):=|x|^A$, where $A\in{\bf R}$, satisfies $f_A\in L^{p,\sigma}(B(0,R))$ 
if and only if 
\begin{equation}
\label{eq:2.3}
 pA+N>0\quad\mbox{for}\quad 1\le\sigma<\infty,\qquad
 pA+N\ge 0\quad\mbox{for}\quad \sigma=\infty.
 \end{equation}
 Furthermore, under condition~\eqref{eq:2.3}, we have 
$$
\|f_A\|_{L^{p,\sigma}(B(0,\sqrt{t}))}\asymp t^{\frac{A}{2}+\frac{N}{2p}}\quad\mbox{for}\quad t\in(0,R^2]
$$
In particular, 
for any $k\in\{0,1,2,\dots\}$, 
by \eqref{eq:1.4} we see that 
\begin{equation}
\label{eq:2.4}
\frac{\|h_k\|_{L^{p,\sigma}(B(0,\sqrt{t}))}}{h_k(\sqrt{t})}\asymp t^{\frac{N}{2p}}\quad\mbox{for}\quad 0<t\le R^2\quad\mbox{if}\quad h_k\in L^{p.\sigma}(B(0,1)).
\end{equation}
%%%%%%%%%%%%
\subsection{Estimates of $h_k$}
%%%%%%%%%%%%
Let $k\in\{0,1,2,\dots\}$. 
Consider the ordinary differential equation
\begin{equation}
\label{eq:2.5}
\frac{d^2}{dr^2}h+\frac{N-1}{r}\frac{d}{dr}h-V_k(r)h=0\quad\mbox{in}\quad(0,\infty).
\end{equation}
Then ODE~\eqref{eq:2.5} has two linearly independent solutions $h_k^+$ and $h_k^-$ such that 
$$
h_k^+(r)=v_{k,\lambda_1}^+(r)(1+o(1)),\qquad
h_k^-(r)=v_{k,\lambda_1}^-(r)(1+o(1)),
$$
as $r\to +0$ and $h_k^-(1)=1$. (See e.g. \cite{IKO}.)
Here
$$
v_{k,\lambda}^+(r):=r^{A^+_{\lambda+\omega_k}},
\quad 
v_{k,\lambda}^-(r):=\left\{
\begin{array}{ll}
r^{-\frac{N-2}{2}}\displaystyle{\left|\log\frac{r}{2}\right|} & \mbox{if $\lambda=\lambda_*$ and $k=0$},\vspace{5pt}\\
r^{A^-_{\lambda+\omega_k}}  & \mbox{otherwise},
\end{array}
\right.
$$
for $\lambda\ge\lambda_*$, where $A^\pm_\lambda$ is as in \eqref{eq:1.3}. 
The solution~$h_k^+$ coincides with $h_k$ (see \eqref{eq:1.4}). 
We first prove the following proposition on the behavior of $h_k$ as $r\to +0$.
\begin{proposition}
\label{Proposition:2.1}
Assume conditions~{\rm ($\mbox{V}$)} and {\rm (N')}. 
Let $R\ge 1$. 
\begin{itemize}
  \item[{\rm (a)}]
  For any $k\in\{0,1,2,\dots\}$ and $\ell\in\{0,1\}$, 
  there exists $C_1>0$ such that 
  $$
  \left|\frac{d^\ell}{dr^\ell}h_k(r)-\frac{d^\ell}{dr^\ell}v^+_{k,\lambda_1}(r)\right|
  \le C_1r^{-\ell+\rho_1}v^+_{k,\lambda_1}(r)
  \quad\mbox{for}\quad r\in(0,R].
  $$
  \item[{\rm (b)}]
  There exists $C_2>0$ and $k_*\in\{0,1,2,\dots\}$ such that 
  $$
  |h_k(r)-v^+_{k,\lambda_1}(r)|\le C_2k^{-1}r^{\rho_1} v^+_{k,\lambda_1}(r),
  \quad
  C_2^{-1}kr^{-1}h_k(r)\le\frac{d}{dr}h_k(r)\le C_2kr^{-1}h_k(r),
  $$
  for $0<r\le R$ and $k\in\{k_*,k_*+1,\dots\}$. 
\end{itemize}
\end{proposition}
For the proof, we prepare the following lemma. 
\begin{lemma}
\label{Lemma:2.1}
Let $k\in\{0,1,2,\dots\}$, $\lambda\ge\lambda_*$ and $R>0$. 
Let $f$ be a continuous function in $(0,R]$ such that 
\begin{equation}
\label{eq:2.6}
|f(r)|\le r^{-2+\epsilon}v^+_{k,\lambda}(r),\qquad r\in(0,R],
\end{equation}
for some $\epsilon>0$. 
Set 
$$
F_{k,\lambda}^+[f](r):=
v^+_{k,\lambda}(r)\int_0^r s^{1-N}[v_{k,\lambda}^+(s)]^{-2}\left(\int_0^s \tau^{N-1}v_{k,\lambda}^+(\tau)f(\tau)\,d\tau\right)\,ds
$$
for $r\in(0,R]$. 
Then there exists $C>0$, independent of $k$, such that 
\begin{equation}
\label{eq:2.7}
\left|\frac{d^\ell}{dr^\ell}F_{k,\lambda}^+[f](r)\right|\le C(k+1)^{\ell-1}r^{-\ell+\epsilon}v^+_{k,\lambda}(r)
\quad\mbox{for}\quad r\in(0,R],
\end{equation}
where $\ell\in\{0,1\}$.
Furthermore, 
\begin{equation}
\label{eq:2.8}
\frac{d^2}{dr^2}F_{k,\lambda}^+[f](r)+\frac{N-1}{r}\frac{d}{dr}F_{k,\lambda}^+[f](r)-\frac{\lambda+\omega_k}{r^2}F_{k,\lambda}^+[f](r)
=f(r)\quad\mbox{in}\quad(0,R]. 
\end{equation}
\end{lemma}
{\bf Proof.}
Assume \eqref{eq:2.6} for some $\epsilon>0$. 
Since $\lambda\ge\lambda_*=-(N-2)^2/4$ and $\omega_k\ge k^2$, 
we see that
$$
N-3+\epsilon+2A^+_{\lambda+\omega_k}=-1+\epsilon+\sqrt{(N-2)^2+4(\lambda+\omega_k)}
\ge -1+\epsilon+2k
$$
for $k\in\{0,1,2,\dots\}$. 
It follows that
\begin{equation*}
\begin{split}
 & s^{1-N}[v_{k,\lambda}^+(s)]^{-2}\int_0^s \tau^{N-1}v_{k,\lambda}^+(\tau)|f(\tau)|\,d\tau\\
 & \le s^{1-N-2A^+_{\lambda+\omega_k}}\int_0^s\tau^{N-3+\epsilon+2A^+_{\lambda+\omega_k}}\,d\tau
\le C(\epsilon+2k)^{-1}s^{-1+\epsilon}\quad\mbox{for}\quad s\in(0,R].
\end{split}
\end{equation*}
Then we easily obtain \eqref{eq:2.7} and \eqref{eq:2.8}, 
and the proof is complete.
$\Box$\vspace{3pt}
\newline
{\bf Proof of Proposition~\ref{Proposition:2.1}.}
Let $V_{\lambda_1}(r):=V(r)-\lambda_1 r^{-2}$ and $R\ge 1$. 
By condition~($\mbox{V}$)~(ii) we find $C_V>0$ such that
\begin{equation}
\label{eq:2.9}
|V_{\lambda_1}(r)|\le C_Vr^{-2+\rho_1}\quad\mbox{for}\quad r>0.
\end{equation}
Define $\{z_n\}_{n=1}^\infty$ inductively by 
$$
z_1(r):=v_{k,\lambda_1}^+(r),
\qquad
z_{n+1}(r):=v_{k,\lambda_1}^+(r)+F_{k,\lambda_1}^+[V_{\lambda_1}z_n](r),\qquad n=1,2,\dots.
$$
We prove that there exists $C_1>0$, independent of $R$ and $k$, such that
\begin{equation}
\label{eq:2.10}
|z_{n+1}(r)-z_n(r)|\le \left(\frac{C_1C_VR^{\rho_1}}{k+1}\right)^nv_{k,\lambda_1}^+(r),\quad r\in(0,R],
\end{equation}
for $n=1,2,\dots$.
By \eqref{eq:2.9} we apply Lemma~\ref{Lemma:2.1} to find $C_*>0$, independent of $k$, such that
$$
|z_2(r)-z_1(r)|\le F_{k,\lambda_1}^+[|V_{\lambda_1}|v_{k,\lambda_1}^+](r)
\le \frac{C_*C_V}{k+1}R^{\rho_1}v_{k,\lambda_1}^+(r),\quad r\in(0,R],
$$
which implies \eqref{eq:2.10} with $n=1$. 
If \eqref{eq:2.10} holds for some $n_*\in\{1,2,\dots\}$, 
then we apply Lemma~\ref{Lemma:2.1} again to obtain 
$$
|z_{n_*+2}(r)-z_{n_*+1}(r)|\le F_{k,\lambda_1}^+[|V_{\lambda_1}||z_{n_*+1}-z_{n_*}|](r)
\le\left(\frac{C_*C_VR^{\rho_1}}{k+1}\right)^{n_*+1}v_{k,\lambda_1}^+(r)
$$
for $r\in(0,R]$. Thus inequality~\eqref{eq:2.10} holds for $n=1,2,\dots$ with $C_1=C_*$.

We prove assertion~(a). Let $k\in\{0,1,2,\dots\}$ and fix it. 
By the regularity of $h_k$, it suffices to treat the case when $R>0$ is small enough. 
By \eqref{eq:2.10}, taking small enough $R>0$, we have 
$$
|z_{n+1}(r)-z_n(r)|\le\left(\frac{1}{2}\right)^nv_{k,\lambda_1}^+(r),\quad r\in(0,R],
$$
for $n=1,2,\dots$.
Then, applying the standard theory for ordinary differential equations, 
we see that the limit function 
$z(r):=\lim_{n\to\infty}z_n(r)$ exists in $(0,R]$ and 
\begin{equation}
\label{eq:2.11}
z(r)=v_{k,\lambda_1}^+(r)+F^+_{k,\lambda_1}[V_{\lambda_1}z](r),
\qquad
|z(r)|\le 2v_{k,\lambda_1}^+(r),
\end{equation}
for $r\in(0,R]$. 
By Lemma~\ref{Lemma:2.1} we find $C_2>0$, independent of $k$, such that 
\begin{equation}
\label{eq:2.12}
\left|\frac{d^\ell}{dr^\ell}F^+_{k,\lambda_1}[V_{\lambda_1}z](r)\right|
\le C_2(k+1)^{\ell-1}r^{-\ell+\rho_1}v_{k,\lambda_1}^+(r)
\end{equation}
for $r\in(0,R]$, where $\ell\in\{0,1\}$. 
Furthermore, we see that $z$ is a solution to ODE~\eqref{eq:2.5} in $(0,R]$. 
Recalling that $h_k^\pm$ are linearly independent solutions to ODE~\eqref{eq:2.5}, 
we find $a$, $b\in{\bf R}$ such that 
$$
z(r)=ah_k^+(r)+bh_k^-(r)\quad\mbox{for}\quad r\in(0,R].
$$
Since $h^+_k(r)/h^-_k(r)\to +0$ as $r\to+0$ and $z(r)=v^+_{k,\lambda_1}(r)(1+o(1))$ as $r\to +0$, 
we see that $a=1$ and $b=0$, that is, $z(r)=h_k^+(r)=h_k(r)$ for $r\in(0,R]$. 
This together with \eqref{eq:2.12} implies assertion~(a) for small enough $R>0$. 
Thus assertion~(a) follows. 

The proof of assertion~(b) is similar. 
Let $R>0$ and fix it. 
By \eqref{eq:2.10}, taking large enough $k_*\in\{0,1,2,\dots\}$, we have 
$$
|z_{n+1}(r)-z_n(r)|\le\left(\frac{1}{2}\right)^nv_{k,\lambda_1}^+(r),\quad r\in(0,R],
$$
for $n=1,2,\dots$ and $k\ge k_*$.
Similarly to assertion~(a), 
we see that the limit function $z(r):=\lim_{n\to\infty}z_n(r)$ exists in $(0,R]$ and 
$z$ satisfies \eqref{eq:2.11}, \eqref{eq:2.12} and $z(r)=h_k(r)$ for $r\in(0,R]$ and $k\ge k_*$.
Furthermore, by \eqref{eq:2.11} and \eqref{eq:2.12} we see that
$$
\left|\frac{d}{dr}h_k(r)-A^+_{\lambda_1+\omega_k}r^{-1}v^+_{k,\lambda_1}(r)\right|
\le C_2r^{-1+\rho_1}v_{k,\lambda_1}^+(r)
$$
for $r\in(0,R]$. Since $A^+_{\lambda_1+\omega_k}=k(1+o(1))$ as $k\to\infty$ (see \eqref{eq:1.3}), 
taking large enough $k_*$ if necessary, 
we obtain 
$$
\frac{d}{dr}h_k(r)\asymp A^+_{\lambda_1+\omega_k}r^{-1}v^+_{k,\lambda_1}(r)
\asymp kr^{-1}v^+_{k,\lambda_1}(r)
$$
for $r\in(0,R]$ and large enough $k\ge k_*$. Then we complete the proof of assertion~(b). 
Thus Proposition~\ref{Proposition:2.1} follows.
$\Box$
\vspace{5pt}

Next we prove the following proposition on the behavior of $h_k$ as $r\to\infty$.
\begin{proposition}
\label{Proposition:2.2}
Assume conditions~{\rm ($\mbox{V}$)} and {\rm (N')}. 
Let $v_k$ be as in \eqref{eq:1.5}.
\begin{itemize}
  \item[{\rm (a)}] 
  Let $k\in\{0,1,2,\dots\}$ and $\ell\in\{0,1\}$. Then
  $$
  \frac{d^\ell}{dr^\ell}h_k(r)
  =c_k(1+o(1))\frac{d^\ell}{dr^\ell}v_k(r)+o(r^{-\ell}v_k(r))
  \quad\mbox{as}\quad r\to\infty.
  $$
  \item[{\rm (b)}] 
  There exist $C>0$ and $k_*\in\{0,1,2,\dots\}$ such that
  $$
  C^{-1}\le\frac{h_k(r)}{v_k(r)}\le C\quad\mbox{in}\quad [1,\infty),
  \quad
  C^{-1}k r^{-1}v_k(r)\le\frac{d}{dr}h_k(r)\le Ck r^{-1}v_k(r), 
  $$
  for $r\ge 1$ and $k\in\{k_*,k_*+1,\dots\}$.
\end{itemize}
\end{proposition}
{\bf Proof.}
Let $V_{\lambda_2}(r):=V(r)-\lambda_2r^{-2}$, $k\in\{0,1,2,\dots\}$ and $\ell\in\{0,1\}$. 
We prove assertion~(a).
\vspace{3pt}
\newline
\underline{Step 1}: 
Assume either 
\begin{equation}
\label{eq:2.13}
k\ge 1\quad\mbox{or}\quad \mbox{$H$ is subcritical and $\lambda>\lambda_*$}. 
\end{equation}
It follows from \eqref{eq:1.5} that
$h_k(r)=c_kv^+_{k,\lambda_2}(r)(1+o(1))$ as $r\to\infty$.
Let $R\ge 1$. 
Let $a_{k,R}$, $b_{k,R}\in{\bf R}$ be such that 
\begin{equation*}
\begin{split}
 & h_k(R)=a_{k,R}v^+_{k,\lambda_2}(R)+b_{k,R}v^-_{k,\lambda_2}(R),\\
 & \frac{d}{dr}h_k(R)=a_{k,R}\frac{d}{dr}v^+_{k,\lambda_2}(R)+b_{k,R}\frac{d}{dr}v^-_{k,\lambda_2}(R).
\end{split}
\end{equation*}
Set 
$z_0(r):=a_{k,R}v^+_{k,\lambda_2}(r)+b_{k,R}v^-_{k,\lambda_2}(r)\quad\mbox{for}\quad r\in(0,\infty)$.
Then 
\begin{equation}
\label{eq:2.14}
|z_0(r)|\le m_{k,R}v^+_{k,\lambda_2}(r)\quad\mbox{for}\quad r\in[1,\infty),
\end{equation}
where $m_{k,R}:=|a_{k,R}|+|b_{k,R}|$.
Define $\{z_n\}_{n=0}^\infty$ inductively by 
$$
z_{n+1}(r):=z_0(r)+G_R[z_n](r)\quad\mbox{for}\quad r\in(0,\infty),
$$
where
$$
G_R[z_n](r):=v^-_{k,\lambda_2}(r)\int_R^r s^{1-N}[v_{k,\lambda_2}^-(s)]^{-2}\left(\int_R^s \tau^{N-1}v_{k,\lambda_2}^-(\tau)
V_{\lambda_2}(\tau)z_n(\tau)\,d\tau\right)\,ds.
$$
On the other hand, by condition~($\mbox{V}$)~(ii) 
we find $C_V'>0$ such that
\begin{equation}
\label{eq:2.15}
|V_{\lambda_2}(r)|\le C_V'r^{-2-\rho_2}\quad\mbox{for}\quad r\in[1,\infty). 
\end{equation}
Then we have
\begin{equation}
\label{eq:2.16}
|z_{n+1}(r)-z_n(r)|\le m_{k,R}\left(\frac{C_V'I_R}{\sqrt{D_k}}\right)^{n+1} v^+_{k,\lambda_2}(r)
\end{equation}
for $r\in[R,\infty)$ and $k\in\{0,1,2,\dots\}$, where 
$$
I_R:=\int_R^\infty \tau^{-1-\rho_2}\,d\tau,
\quad D_k:=D_{\lambda_2+\omega_k}=(N-2)^2+4(\lambda_2+\omega_k)\ge 4k^2.
$$
Indeed, it follows from \eqref{eq:2.15} that 
$$
\int_R^s\tau^{N-1}v_{k,\lambda_2}^-(\tau)
\left|V_{\lambda_2}(\tau)\right|v^+_{k,\lambda_2}(\tau)\,d\tau
\le C_V'\int_R^s\tau^{-1-\rho_2}\,d\tau\le C_V'I_R
$$
for $s\in[R,\infty)$. This implies that 
\begin{equation}
\label{eq:2.17}
\begin{split}
 |G_R[v^+_{k,\lambda_2}](r)|
  & \le C_V'I_Rv^-_{k,\lambda_2}(r)\int_R^r s^{1-N}[v_{k,\lambda_2}^-(s)]^{-2}\,ds\\
 & =C_V'I_Rv^-_{k,\lambda_2}(r)\int_R^r s^{-1+\sqrt{D_k}}\,ds
 \le\frac{C_V'I_R}{\sqrt{D_k}}v^+_{k,\lambda_2}(r)
\end{split}
\end{equation}
for $r\in[R,\infty)$. This together with \eqref{eq:2.14} implies \eqref{eq:2.16} with $n=0$. 
Repeating this argument, 
by induction we see that \eqref{eq:2.16} holds for $n=0,1,2,\dots$. 

Since $I_R\to 0$ as $R\to\infty$, 
we take large enough $R_k\ge 1$ so that 
\begin{equation}
\label{eq:2.18}
\frac{C_V'I_{R}}{\sqrt{D_k}}\le\frac{1}{2}\quad\mbox{for}\quad R\ge R_k.
\end{equation}
Then, for any $R\ge R_k$, 
applying the standard theory for ordinary differential equations, 
we see that $z(r):=\lim_{n\to\infty}z_n$ exists for $r\in[R,\infty)$ and $z$ satisfies
\begin{equation}
\label{eq:2.19}
z(r)=z_0(r)+G_R[z](r),\qquad
|z(r)|\le 2m_{k,R_k}v^+_{k,\lambda_2}(r),
\end{equation}
for $r\in[R,\infty)$. 
Furthermore, 
\begin{equation}
\label{eq:2.20}
\begin{split}
 & \frac{d^2}{dr^2}z+\frac{N-1}{r}\frac{d}{dr}z-(\lambda_2+\omega_k)r^{-2}z=V_{\lambda_2}(r)z\quad\mbox{in}\quad[R,\infty),\\
 & z(R)=z_0(R)=h_k(R),\quad \frac{d}{dr}z(R)=\frac{d}{dr}z_0(R)=\frac{d}{dr}h_k(R). 
\end{split}
\end{equation}
On the other hand, 
by \eqref{eq:1.4} we see that 
$h_k$ satisfies relation~\eqref{eq:2.20} with $z$ replaced by $h_k$. 
These imply that $h_k(r)=z(r)$ for $r\in[R,\infty)$. 
Therefore we obtain 
\begin{equation}
\label{eq:2.21}
h_k(r)=a_{k,R}v^+_{k,\lambda_2}(r)+b_{k,R}v^-_{k,\lambda_2}(r)
+G_R[h_k](r)
\end{equation}
for $r\ge R$, where $R\ge R_k$. 
On the other hand, 
under assumption~\eqref{eq:2.13}, 
by \eqref{eq:1.5} we have
$$
|h_k(r)|\le Cv^+_{k,\lambda_2}(r)\quad\mbox{for}\quad r\in[1,\infty).
$$
Then, similarly to \eqref{eq:2.17}, we see that
\begin{equation}
\label{eq:2.22}
\left|\frac{d^\ell}{dr^\ell}G_R[h_k](r)\right|\le CI_Rr^{-\ell}v^+_{k,\lambda_2}(r)
\end{equation}
for $r\in[R,\infty)$. 
Since $h_k(r)=c_kv^+_{k,\lambda_2}(r)(1+o(1))$ as $r\to\infty$ (see \eqref{eq:1.5}) and $I_R\to 0$ as $R\to\infty$,  
by \eqref{eq:2.21} and \eqref{eq:2.22} we see that 
$a_{k,R}\to c_k$ as $R\to\infty$ and obtain
$$
\frac{d^\ell}{dr^\ell}h_k(r)=c_k(1+o(1))\frac{d^\ell}{dr^\ell}v^+_{k,\lambda_2}(r)
+o\left(r^{-\ell}v^+_{k,\lambda_2}(r)\right)
\quad\mbox{as}\quad r\to\infty. 
$$
Thus assertion~(a) follows under assumption~\eqref{eq:2.13}. 
\vspace{3pt}
\newline
\underline{Step 2}: 
Let $k=0$ and assume either
\begin{equation}
\label{eq:2.23}
{\rm (i)}\quad\mbox{$H$ is subcritical and $\lambda_2=\lambda_*$}
\qquad\mbox{or}\qquad {\rm (ii)}\quad\mbox{$H$ is critical}.
\end{equation} 
Set 
$$
G(r):=v^+_{0,\lambda_2}(r)\int^\infty_r s^{1-N}[v_{0,\lambda_2}^+(s)]^{-2}\left(\int^\infty_s \tau^{N-1}v_{0,\lambda_2}^+(\tau)
V_{\lambda_2}(\tau)h_0(\tau)\,d\tau\right)\,ds.
$$
It follows from \eqref{eq:1.5} and \eqref{eq:2.15} that 
\begin{equation}
\label{eq:2.24}
\tau^{N-1}v_{0,\lambda_2}^+(\tau)
\left|V_{\lambda_2}(\tau)\right|h_0(\tau)
\le C\tau^{N-1}\tau^{A^+_{\lambda_2}}\tau^{-2-\rho_2}\tau^{A^-_{\lambda_2}}(\log \tau)^{B_0}
=C\tau^{-1-\rho_2}(\log \tau)^{B_0}
\end{equation}
for $\tau\ge 2$.
Then we have 
\begin{equation}
\label{eq:2.25}
|G(r)|\le Cv^+_{0,\lambda_2}(r)\int_r^\infty s^{-1-\sqrt{D_0}}s^{-\rho_2}(\log r)^{B_0}\,ds
 \le Cr^{-\rho_2}v_{0,\lambda_2}^-(r)
\end{equation}
for $r\ge 2$. This implies that $G(r)=o(h_0(r))$ as $r\to\infty$. 
Furthermore, $G$ satisfies
$$
\frac{d^2}{dr^2}G+\frac{N-1}{r}\frac{d}{dr}G-\lambda_2r^{-2}G
=V_{\lambda_2}(r)h_0(r)\quad\mbox{in}\quad(0,\infty).
$$
Therefore, setting $\tilde{h}:=h_0-G$, we have 
\begin{equation}
\label{eq:2.26}
\frac{d^2}{dr^2}\tilde{h}+\frac{N-1}{r}\frac{d}{dr}\tilde{h}-\lambda_2r^{-2}\tilde{h}=0\quad\mbox{in}\quad(0,\infty).
\end{equation}
Since $v^\pm_{0,\lambda_2}$ are linearly independent solutions to ODE~\eqref{eq:2.26}, 
we find $a$, $b\in{\bf R}$ such that 
$$
h_0(r)-G(r)=\tilde{h}(r)=av^+_{0,\lambda_2}(r)+bv^-_{0,\lambda_2}(r)\quad\mbox{for}\quad r>0. 
$$
On the other hand, it follows from \eqref{eq:1.5} 
that $h_0(r)=c_0v^-_{0,\lambda_2}(r)(1+o(1))$ as $r\to\infty$. 
This implies that $a=0$ and $b=c_0$, that is, 
$$
h_0(r)=c_0v^-_{0,\lambda_2}(r)+G(r)\quad\mbox{in}\quad(0,\infty).
$$
Then assertion~(a) easily follows from \eqref{eq:2.24} and \eqref{eq:2.25} under assumption~\eqref{eq:2.23}. 
The proof of assertion~(a) is complete. 
\vspace{3pt}
\newline
\underline{Step 3}: 
We prove assertion~(b). 
In this step the letter $C$ denotes generic positive constants independent of $k$. 
Let $R=1$. 
Since $D_k\to\infty$ as $k\to\infty$, 
taking large enough $k_*\in\{1,2,\dots\}$, 
we have 
$$
\frac{C_V'I_1}{\sqrt{D_k}}\le Ck^{-1}\le\frac{1}{2}\quad\mbox{for}\quad k\ge k_*,
$$
instead of \eqref{eq:2.18}.
Then, similarly to \eqref{eq:2.19} (with $R_k=1$), 
we see that $z(r):=\lim_{n\to\infty}z_n$ exists for $r\in[1,\infty)$ and $z$ satisfies
$$
z(r)=z_0(r)+G_1[z](r),\qquad
|z(r)|\le 2m_{k,1}v^+_{k,\lambda_2}(r), 
$$
for $r\in[1,\infty)$, where $k\ge k_*$. Furthermore, $z=h_k$ in $[1,\infty)$ and
\begin{equation*}
\begin{split}
 & h_k(r)=a_{k,1}v^+_{k,\lambda_2}(r)+b_{k,1}v^-_{k,\lambda_2}(r)
+G_1[h_k](r),\\
 & |G_1[h_k](r)|\le Ck^{-1}v^+_{k,\lambda_2}(r),\quad
 \left|\frac{d}{dr}G_1[h_k](r)\right|\le Cr^{-1}v^+_{k,\lambda_2}(r),
\end{split}
\end{equation*}
for $r\in[1,\infty)$. 
On the other hand, 
taking large enough $k_*$ if necessary, 
by Proposition~\ref{Proposition:2.1} 
we see that
$$
C^{-1}<h_k(1)\le C,\qquad C^{-1}k\le \frac{d}{dr}h_k(1)\le Ck,
$$
for $k\ge k_*$. These imply that
\begin{equation}
\label{eq:2.27}
C^{-1}<a_{k,1}+b_{k,1}\le C,\quad C^{-1}k\le A^+_{\lambda_2+\omega_k}a_{k,1}+A^-_{\lambda_2+\omega_k}b_{k,1}\le Ck,
\end{equation}
for $k\ge k_*$. 
Since 
$$
D_k=(N-2)^2+4(\lambda_2+\omega_k)=4k^2(1+o(1))\quad\mbox{as}\quad k\to\infty,
$$
it follows from \eqref{eq:2.27} that
\begin{equation*}
\begin{split}
k & \asymp A^+_{\lambda_2+\omega_k}a_{k,1}+A^-_{\lambda_2+\omega_k}b_{k,1}
+\frac{\sqrt{D_k}}{2}(a_{k,1}+b_{k,1})\\
 & =(A^+_{\lambda_2+\omega_k}-A^-_{\lambda_2+\omega_k})a_{k,1}+A^-_{\lambda_2+\omega_k}(a_{k,1}+b_{k,1})
+\frac{\sqrt{D_k}}{2}(a_{k,1}+b_{k,1})\\
 & =\sqrt{D_k}\,a_{k,1}-\frac{N-2}{2}(a_{k,1}+b_{k,1})
 =2k(1+o(1))a_{k,1}+O(1)\quad\mbox{as}\quad k\to\infty.
\end{split}
\end{equation*}
This implies that 
$$
C^{-1}\le a_{k,1}\le C
$$
for large enough $k$. 
Then, by \eqref{eq:2.27} we see that $|b_{k,1}|\le C$ for large enough~$k$.
Therefore we find $R\ge1$ such that 
\begin{equation}
\label{eq:2.28}
h_k(r)\asymp v_{k,\lambda_2}^+(r)=v_k(r)
\end{equation}
for $r\ge R$ and large enough $k$. Furthermore
\begin{equation*}
\begin{split}
 & C^{-1}kr^{-1}v^+_{k,\lambda_2}(r)-Ckr^{-1}v^-_{k,\lambda_2}(r)\\
 & \le \frac{d}{dr}h_k(r)=A^+_{\lambda_2+\omega_k}a_{k,1}r^{-1}v^+_{k,\lambda_2}(r)+A^-_{\lambda_2+\omega_k}b_{k,1}r^{-1}v^-_{k,\lambda_2}(r)
+\frac{d}{dr}G_1[h_k](r)\\
 & \le Ckr^{-1}v^+_{k,\lambda_2}(r)+Ckr^{-1}v^-_{k,\lambda_2}(r)
\end{split}
\end{equation*}
for $r\ge 1$ and large enough $k$. 
Taking large enough $R$ if necessary, we see that
\begin{equation}
\label{eq:2.29}
\frac{1}{2}C^{-1}kr^{-1}v^+_{k,\lambda_2}(r)\le \frac{d}{dr}h_k(r)\le 2Ckr^{-1}v^+_{k,\lambda_2}(r)
\end{equation}
for $r\in[R,\infty)$ and large enough $k$.
Combining Proposition~\ref{Proposition:2.1}~(b) with \eqref{eq:2.28} and \eqref{eq:2.29},
we complete the proof of assertion~(b).
Thus Proposition~\ref{Proposition:2.2} follows. 
$\Box$\vspace{5pt}

Combining Propositions~\ref{Proposition:2.1} and \ref{Proposition:2.2},
we have:
\begin{proposition}
\label{Proposition:2.3} 
Assume conditions~{\rm (V)} and {\rm (N')}. 
Then there exists $C>0$ such that
$$
C^{-1}\le\displaystyle{\frac{h_k(r)}{v_{k,\lambda_1}^+(r)}}\le C
\quad\mbox{in}\quad(0,1],\qquad
C^{-1}\le\displaystyle{\frac{h_k(r)}{v_k(r)}}\le C\quad\mbox{in}\quad(1,\infty),
$$
for $k\in\{0,1,2,\dots\}$. 
\end{proposition}

At the end of this subsection we prove the following proposition.
\begin{proposition}
\label{Proposition:2.4} 
Assume conditions~{\rm (V)} and {\rm (N')}.
\begin{itemize}
  \item[{\rm (a)}] 
  Let $(p,p,\sigma,\sigma)\in\Lambda$ be such that $h_0\in L^{p,\sigma}(B(0,1))$. 
  Then there exists $C_1>0$ such that 
  $$
  \frac{\|h_0\|_{L^{p,\sigma}(B(0,\sqrt{t}))}}{h_0(\sqrt{t})}\ge C_1t^{\frac{N}{2p}}\quad\mbox{for}\quad t>0.
  $$ 
  \item[{\rm (b)}] 
  There exists $C_2>0$ such that 
  $$
  \int_0^r s^{N-1}h_k(s)^2\,ds\le C_2(k+1)^{-1}r^N h_k(r)^2
  $$
  for $r>0$ and $k\in\{0,1,2,\dots\}$. 
\end{itemize}
\end{proposition}
{\bf Proof.} 
By \eqref{eq:1.4} and \eqref{eq:1.5} we find $C>0$ and $\epsilon\in(0,1)$ such that 
$$
\frac{\|h_0^\epsilon\|_{L^1(B(0,\sqrt{t}))}}{h_0(\sqrt{t})^\epsilon}\ge Ct^{\frac{N}{2}}\quad\mbox{for}\quad t>0. 
$$
It follows from \eqref{eq:2.1} that 
$$
\|h_0^\epsilon\|_{L^1(B(0,\sqrt{t}))}\le |B(0,\sqrt{t})|^{1-\frac{\epsilon}{p}}\|h_0\|_{L^{p,\sigma}(B(0,\sqrt{t}))}^\epsilon
\quad\mbox{for}\quad t>0. 
$$
These imply assertion~(a). 
Assertion~(b) follows from Proposition~\ref{Proposition:2.3} 
(see also \cite[(3.7)]{IM}).  
$\Box$
%%%%%%%%%%%%%%%%%%%%%%%%%%%%%%%%%%%%
%%%%%%%%%%%%%%%%%%%%%%%%%%%%%%%%%%%%
\section{$L^{p,\sigma}$-$L^{q,\theta}$ estimates of $e^{-tH}$ }
%%%%%%%%%%%%%%%%%%%%%%%%%%%%%%%%%%%%
%%%%%%%%%%%%%%%%%%%%%%%%%%%%%%%%%%%%
The aim of this section is to prove the following proposition.
\begin{proposition}
\label{Proposition:3.1}
Assume conditions~{\rm (V)} and {\rm (N')}. Let $(p,q,\sigma,\theta)\in\Lambda$. 
\begin{itemize}
  \item[{\rm (a)}] 
  Let $\ell\in\{0,1\}$, $j\in\{0,1,2,\dots\}$ and $\delta\in (0,1]$.  
  Then there exist $C_1>0$ and $C_2>0$ such that 
  \begin{equation*}
  \begin{split}
  \left\|\partial_t^j\nabla^\ell e^{-tH}\phi\right\|_{L^{q,\theta}(B(0,\delta\sqrt{t})^c)}
   & \le C_1t^{-\frac{N}{2}\left(1-\frac{1}{q}\right)-\frac{\ell}{2}-j}
  \left[\frac{\|h_0\phi\|_{L^1(B(0,\sqrt{t}))}}{h_0(\sqrt{t})}+t^{\frac{N}{2p'}}\|\phi\|_{L^{p,\sigma}(B(0,\delta\sqrt{t})^c)}\right]\\
   & \le C_2t^{-\frac{N}{2}\left(1-\frac{1}{q}\right)-\frac{\ell}{2}-j}
   \frac{\|h_0\|_{L^{p',\sigma'}(B(0,\sqrt{t}))}}{h_0(\sqrt{t})}\|\phi\|_{L^{p,\sigma}}
  \end{split}
  \end{equation*}
  for $\phi\in C_c({\bf R}^N)$ and $t>0$. 
  \item[{\rm (b)}] 
  There exist $C_3>0$ and $C_4>0$ such that 
  \begin{equation*}
  \begin{split}
  \|e^{-tH}\phi\|_{L^{q,\theta}}
   & \le C_3t^{-\frac{N}{2}}\frac{\|h_0\|_{L^{q,\theta}(B(0,\sqrt{t}))}}{h_0(\sqrt{t})}
  \left[\frac{\|h_0\phi\|_{L^1(B(0,\sqrt{t}))}}{h_0(\sqrt{t})}
  +t^{\frac{N}{2p'}}\|\phi\|_{L^{p,\sigma}(B(0,\sqrt{t})^c)}\right]\\
   & \le C_4t^{-\frac{N}{2}}\frac{\|h_0\|_{L^{q,\theta}(B(0,\sqrt{t}))}\|h_0\|_{L^{p',\sigma'}(B(0,\sqrt{t}))}}{h_0(\sqrt{t})^2}\|\phi\|_{L^{p,\sigma}}
  \end{split}
  \end{equation*}
  for $\phi\in C_c({\bf R}^N)$ and $t>0$. 
\end{itemize}
\end{proposition}

We first recall the following lemma 
on an upper Gaussian estimate of $p=p(x,y,t)$. 
See \cite{IKO}*{Theorem~1.3}. 
\begin{lemma}
\label{Lemma:3.1}
Assume conditions {\rm (V)} and~{\rm (N')}. 
Then there exists $C>0$ such that
$$
0<p(x,y,t)\le Ct^{-\frac{N}{2}}\frac{\tilde{h}_0(x,t)\tilde{h}_0(y,t)}{h_0(\sqrt{t})^2}\exp\left(-\frac{|x-y|^2}{Ct}\right)
$$
for $x$, $y\in{\bf R}^N\setminus\{0\}$ and $t>0$, where 
$\tilde{h}_0(x,t):=h_0(\min\{|x|,\sqrt{t}\})$. 
\end{lemma}
Combining Lemma~\ref{Lemma:3.1} and the parabolic regularity theorems, 
we have: 
\begin{lemma}
\label{Lemma:3.2}
Assume conditions~{\rm (V)} and {\rm (N')}. 
Let $\ell\in\{0,1\}$, $j\in\{0,1,2,\dots\}$ and $\delta\in(0,1]$. 
Then there exists $C>0$ such that
\begin{equation*}
|\partial_t^j \nabla^\ell p(x,y,t)|
\le Ct^{-\frac{N}{2}-\frac{\ell}{2}-j}
\frac{\tilde{h}_0(x,t)\tilde{h}_0(y,t)}{h_0(\sqrt{t})^2}\exp\left(-\frac{|x-y|^2}{Ct}\right)
\end{equation*}
for $x\in B(0,\delta\sqrt{t})^c$, $y\in \R^N\setminus\{0\}$ and $t>0$.
\end{lemma}
{\bf Proof.}
Let $y\in{\bf R}^N\setminus\{0\}$, $t>0$ and $\delta\in(0,1]$. 
Set $k=\delta\sqrt{t}/2$. 
Let $x\in B(0,\delta\sqrt{t})^c$ and  set $p_k(z,s):=p(x+kz,y,t+k^2s)$ for $(z,s)\in B(0,1)\times(-1,1)$. 
Then $p_k$ satisfies
\begin{equation}
\label{eq:3.1}
\partial_sp_k-\Delta_zp_k+\tilde{V}(z)p_k=0 \quad\mbox{in}\quad (x,t) \in B(0,1)\times(-1,1),
\end{equation}
where $\tilde{V}(z)=k^2V(|x+kz|)$. 
By condition~(V)~(iii) we have 
\begin{equation}
\label{eq:3.2}
|\nabla\tilde{V}(z)|\le Ck^{3}|x+kz|^{-3}\le C
\quad\mbox{for}\quad z\in B(0,1).
\end{equation}
Let $\ell\in\{0,1\}$ and $j\in\{0,1,2,\dots\}$. 
By \eqref{eq:3.2} 
we apply the parabolic regularity theorems to \eqref{eq:3.1} 
and obtain
$$
k^{\ell+2j}|\partial_t^j\nabla^\ell p(x,y,t)|
=|\partial_t^j\nabla^\ell p_k(0,0)|
\le C\|p_k\|_{L^\infty(B(0,1)\times(-1,1))}.
$$
Then, due to the relation that $h_0(r/2)\asymp h_0(r)\asymp h_0(2r)$ for $r>0$ (see \eqref{eq:1.4} and \eqref{eq:1.5}),
by Lemma~\ref{Lemma:3.1} we obtain the desired inequality.
$\Box$ 
\vspace{5pt}%

We prove Proposition~\ref{Proposition:3.1} by using Lemmas~\ref{Lemma:3.1} and \ref{Lemma:3.2}. 
\vspace{5pt}
\newline
{\bf Proof of Proposition~\ref{Proposition:3.1}.}
Let $(p,q,\sigma,\theta)\in\Lambda$ and $\phi\in L^{p,\sigma}$. 
By Lemmas~\ref{Lemma:3.1} and \ref{Lemma:3.2} we find $c>0$ such that 
\begin{equation}
\label{eq:3.3}
   \left|\left[\partial_t^j\nabla^\ell e^{-tH}\phi\right](x)\right|
   \le Ct^{-\frac{\ell}{2}-j}[I(x,t)+J(x,t)],\qquad t>0,
\end{equation}
for $x\in B(0,\delta\sqrt{t})^c$ if $(j,\ell)\not=(0,0)$ 
and $x\in{\bf R}^N\setminus\{0\}$ if $(j,\ell)=(0,0)$, 
where
\begin{equation*}
\begin{split}
 & I(x,t):=h_0(\sqrt{t})^{-2}\tilde{h}_0(x,t)\int_{B(0,\sqrt{t})}h_0(|y|)G_c(x-y,t)|\phi(y)|\,dy,\\
 & J(x,t):=h_0(\sqrt{t})^{-1}\tilde{h}_0(x,t)\int_{B(0,\sqrt{t})^c}G_c(x-y,t)|\phi(y)|\,dy,\\
 & G_c(x,t):=t^{-\frac{N}{2}}\exp\left(-\frac{|x|^2}{ct}\right).
\end{split}
\end{equation*}
We prove assertion~(a). 
Let $\delta\in(0,1]$, $x\in B(0,\delta\sqrt{t})^c$ and $t>0$. 
Thanks to \eqref{eq:1.4} and \eqref{eq:1.5}, 
we have $\tilde{h}_0(x,t)\le Ch_0(\sqrt{t})$. 
Then it follows from \eqref{eq:2.2} that 
\begin{equation}
\label{eq:3.4}
\begin{split}
\|I(t)\|_{L^{q,\theta}(B(0,\delta\sqrt{t})^c)}
 & \le Ch_0(\sqrt{t})^{-1}
\left\|\int_{{\bf R}^N}G_c(\cdot\,-y,t)h_0(|y|)|\phi(y)|\chi_{B(0,\sqrt{t})}(y)\,dy\right\|_{L^{q,\theta}}\\
 & \le Ch_0(\sqrt{t})^{-1}\|G_c(t)\|_{L^{q,\theta}}\|h_0\phi\|_{L^1(B(0,\sqrt{t}))}\\
 & \le Ct^{-\frac{N}{2}\left(1-\frac{1}{q}\right)}h_0(\sqrt{t})^{-1}\|h_0\phi\|_{L^1(B(0,\sqrt{t}))}.
\end{split}
\end{equation}
Similarly, by \eqref{eq:2.2}
we obtain
\begin{equation}
\label{eq:3.5}
\begin{split}
\|J(t)\|_{L^{q,\theta}(B(0,\delta\sqrt{t})^c)}
 & \le\|G_c(t)\|_{L^{r,s}}\|\phi\|_{L^{p,\sigma}(B(0,\delta\sqrt{t})^c)}\\
 & \le Ct^{-\frac{N}{2}\left(1-\frac{1}{r}\right)}\|\phi\|_{L^{p,\sigma}(B(0,\delta\sqrt{t})^c)}
 =Ct^{-\frac{N}{2}\left(\frac{1}{p}-\frac{1}{q}\right)}\|\phi\|_{L^{p,\sigma}(B(0,\delta\sqrt{t})^c)}
\end{split}
\end{equation}
for $t>0$, where $1\le r\le\infty$ and $1\le s\le\infty$ with 
$$
\frac{1}{r}+\frac{1}{p}=\frac{1}{q}+1,\qquad \frac{1}{\theta}=\frac{1}{s}+\frac{1}{\sigma}.
$$
Combining \eqref{eq:3.3}, \eqref{eq:3.4} and \eqref{eq:3.5}, 
we see that
$$
\|\partial_t^j\nabla^\ell e^{-tH}\phi\|_{L^{q,\theta}(B(0,\delta\sqrt{t})^c)}
\le Ct^{-\frac{N}{2}\left(1-\frac{1}{q}\right)-\frac{\ell}{2}-j}
\left[\frac{\|h_0\phi\|_{L^1(B(0,\sqrt{t}))}}{h_0(\sqrt{t})}+t^{\frac{N}{2p'}}\|\phi\|_{L^{p,\sigma}(B(0,\delta\sqrt{t})^c)}\right],
$$
which together with Proposition~\ref{Proposition:2.4}~(a) and \eqref{eq:2.1} implies assertion~(a).

We prove assertion~(b). 
Let $x\in B(0,\sqrt{t})$ and $t>0$. 
It follows from $\tilde{h}_0(x,t)=h_0(|x|)$ that
$$
I(x,t)\le h_0(\sqrt{t})^{-2}h_0(|x|)\int_{B(0,\sqrt{t})}
h_0(|y|)G_c(x-y,t)|\phi(y)|\,dy,
$$
which implies that  
\begin{equation}
\label{eq:3.6}
\begin{split}
\|I(t)\|_{L^{q,\theta}(B(0,\sqrt{t}))}
 & \le C\frac{\|h_0\|_{L^{q,\theta}(B(0,\sqrt{t}))}}{h_0(\sqrt{t})^2}
 \sup_{x\in B(0,\sqrt{t})}\int_{B(0,\sqrt{t})} h_0(|y|)G_c(x-y,t)|\phi(y)|\,dy\\
 & \le C\frac{\|h_0\|_{L^{q,\theta}(B(0,\sqrt{t}))}}{h_0(\sqrt{t})^2}\cdot Ct^{-\frac{N}{2}}\|h_0\phi\|_{L^1(B(0,\sqrt{t}))}
\end{split}
\end{equation}
for $t>0$.
Similarly, we obtain
\begin{equation}
\label{eq:3.7}
\begin{split}
 \|J(t)\|_{L^{q,\theta}(B(0,\sqrt{t}))}
 & \le Ch_0(\sqrt{t})^{-1}\|h_0\|_{L^{q,\theta}(B(0,\sqrt{t}))}\sup_{x\in B(0,\sqrt{t})}\int_{B(0,\sqrt{t})^c}G_c(x-y,t)|\phi(y)|\,dy\\
 & \le Ct^{-\frac{N}{2p}}\frac{\|h_0\|_{L^{q,\theta}(B(0,\sqrt{t}))}}{h_0(\sqrt{t})}\|\phi\|_{L^{p,\sigma}(B(0,\sqrt{t})^c)}
\end{split}
\end{equation}
for $t>0$.
By assertion~(a) with $(j,\ell)=(0,0)$, \eqref{eq:3.6} and \eqref{eq:3.7} 
we have 
$$
\|e^{-tH}\phi\|_{L^{q,\theta}}
\le Ct^{-\frac{N}{2}}
\left[\frac{\|h_0\|_{L^{q,\theta}(B(0,\sqrt{t}))}}{h_0(\sqrt{t})}+t^{\frac{N}{2q}}\right]
\left[\frac{\|h_0\phi\|_{L^1(B(0,\sqrt{t}))}}{h_0(\sqrt{t})}
+t^{\frac{N}{2p'}}\|\phi\|_{L^{p,\sigma}(B(0,\sqrt{t})^c)}\right]
$$
for $t>0$. This together with Proposition~\ref{Proposition:2.4}~(a) and \eqref{eq:2.1} 
implies assertion~(b). 
Thus Proposition~\ref{Proposition:3.1} follows.
$\Box$
%%%%%%%%%%%%%%%%%%%%%%%%%%%%%%%%%%%%
%%%%%%%%%%%%%%%%%%%%%%%%%%%%%%%%%%%%
\section{Radially symmetric solutions}
%%%%%%%%%%%%%%%%%%%%%%%%%%%%%%%%%%%%
%%%%%%%%%%%%%%%%%%%%%%%%%%%%%%%%%%%%
In this section we obtain uniform estimates of time derivatives and spatial derivatives of $v_{k,i}$ 
inside parabolic cones with respect to $k$ and $i$.
%%%%%%%%%%%%%%%%%%%%%%%%%%%
\subsection{Estimates of time derivatives}
%%%%%%%%%%%%%%%%%%%%%%%%%%%
Let $j\in\{0,1,2,\dots\}$.
In this subsection we prove the following proposition on uniform estimates of $\partial_t^jv_{k,i}$.
Proposition~\ref{Proposition:4.1} is one of the main ingredients in the proof of Theorem~\ref{Theorem:1.1}.
\begin{proposition}
\label{Proposition:4.1}
Assume conditions~{\rm ($\mbox{V}$)} and {\rm (N')}. 
Let $(p,p,\sigma,\sigma)\in\Lambda$ be such that $h_0\in L^{p,\sigma}(B(0,1))\cap L^{p',\sigma'}(B(0,1))$. 
Let $v_{k,i}$ be as in \eqref{eq:1.9}. 
Then, for any $j\in\{0,1,2,\dots\}$, 
there exist $C>0$ and $\delta\in(0,1]$ such that 
\begin{equation}
\label{eq:4.1}
|\partial_t^j v_{k,i}(x,t)|\le CM_{k,i}t^{-\frac{N}{2}-j}\frac{\|h_0\|_{L^{p',\sigma'}(B(0,\sqrt{t}))}}{h_k(\delta\sqrt{t})h_0(\sqrt{t})}h_k(|x|)\|\phi\|_{L^{p,\sigma}}
\end{equation}
for $x\in B(0,\delta\sqrt{t})$, $t>0$, $k\in\{0,1,2,\dots\}$ and $i\in\{1,\dots,d_k\}$. 
Here $M_{k,i}:=\|Q_{k,i}\|_{L^\infty({\bf S}^{N-1})}$.
\end{proposition}
For the proof of Proposition~\ref{Proposition:4.1}, 
we start by proving the following lemma.
\begin{lemma}
\label{Lemma:4.1}
Assume the same conditions as in Proposition~{\rm\ref{Proposition:4.1}}. 
Then, for any $\delta\in(0,1]$ and $j\in\{0,1,2,\dots\}$, 
there exists $C>0$ such that
\begin{equation}
\label{eq:4.2}
|\partial_t^j v_{k,i}(x,t+t_0)|\le CM_{k,i}t^{-\frac{N}{2}-j}\frac{\|h_0\|_{L^{p',\sigma'}(B(0,\sqrt{t}))}}{h_0(\sqrt{t})}\|\phi\|_{L^{p,\sigma}}
\end{equation}
for $x\in\partial B(0,\delta\sqrt{t})$, $t>0$, $t_0\in[0,1]$, $k\in\{0,1,2,\dots\}$ and $i\in\{1,\dots,d_k\}$.
\end{lemma}
{\bf Proof.}
Let $\phi\in C_c({\bf R}^N)$, $j\in\{0,1,2,\dots\}$ and $\delta\in(0,1]$.
Let $u(t):=e^{-tH}\phi$ for $t>0$.
By Proposition~\ref{Proposition:3.1}~(a) with~$q=\theta=\infty$, we have
\begin{equation*}
|\partial_t^j u(x,t+t_0)| 
\le Ct^{-\frac{N}{2}-j}
\frac{\|h_0\|_{L^{p',\sigma'}(B(0,\sqrt{t}))}}{h_0(\sqrt{t})}\|u(t_0)\|_{L^{p,\sigma}}
\end{equation*}
for $x\in\partial B(0,\delta\sqrt{t})$, $t>0$ and $t_0\in[0,1]$.
Since $h_0 \in L^{p,\sigma}(B(0,1))\cap L^{p',\sigma'}(B(0,1))$,
by Proposition~\ref{Proposition:3.1}~(b) and \eqref{eq:2.4} 
we have
\begin{equation*}
\|u(t_0)\|_{L^{p,\sigma}}
\le Ct_0^{-\frac{N}{2}}\frac{\|h_0\|_{L^{p,\sigma}(B(0,\sqrt{t_0}))}\|h_0\|_{L^{p',\sigma'}(B(0,\sqrt{t_0}))}}{h_0(\sqrt{t_0})^2}\|\phi\|_{L^{p,\sigma}}
\le C\|\phi\|_{L^{p,\sigma}}
\end{equation*}
for $t_0\in(0,1]$.
These imply that
\begin{equation}
\label{eq:4.3}
|\partial_t^j u(x,t+t_0)|
\le Ct^{-\frac{N}{2}-j}
\frac{\|h_0\|_{L^{p',\sigma'}(B(0,\sqrt{t}))}}{h_0(\sqrt{t})}\|\phi\|_{L^{p,\sigma}}
\end{equation}
for $x\in\partial B(0,\delta\sqrt{t})$, $t>0$ and $t_0\in[0,1]$.
On the other hand, 
it follows from \eqref{eq:1.10} that 
$$
\int_{\partial B(0,1)}\partial_t^ju(|x|y,t)Q_{k,i}(y)\,d\sigma_y=\partial_t^jv_{k,i}(|x|,t)
$$
for $x\in{\bf R}^N\setminus\{0\}$ and $t>0$. 
This together with \eqref{eq:4.3} implies \eqref{eq:4.2}. 
Thus Lemma~\ref{Lemma:4.1} follows. 
$\Box$\vspace{5pt}

Next we prepare the following lemma on supersolutions. 
\begin{lemma}
\label{Lemma:4.2}
Assume the same conditions as in Proposition~{\rm\ref{Proposition:4.1}}. 
Let $k\in\{0,1,2,\dots\}$, $T\ge 0$ and $M>0$. 
Let $\zeta\in C^1((T,\infty))$ be such that 
$$
0\le-\zeta'(t)\le (k+1)M t^{-1}\zeta(t)\quad\mbox{in}\quad(T,\infty).
$$
Set 
\begin{equation*}
\begin{split}
f_k(|x|) & :=\int_0^{|x|} s^{1-N}\nu_k(s)^{-1}\left(\int_0^s \tau^{N-1}\nu_k(\tau)\,d\tau\right)\,ds,\qquad \nu_k(|x|)=h_k(|x|)^2,\\
z_*(x,t) & :=2\zeta(t)\left[1-(k+1)Mt^{-1}f_k(|x|)\right],
\end{split}
\end{equation*}
for $(x,t)\in{\bf R}^N\times(T,\infty)$. 
Then $z_*$ satisfies
$$
\partial_t z_*\ge\frac{1}{\nu_k}\mbox{div}\,(\nu_k\nabla z_*),
\qquad z_*\le 2\zeta(t),
$$
for $(x,t)\in{\bf R}^N\times(T,\infty)$. 
Furthermore, there exists $\delta>0$ such that 
\begin{equation}
\label{eq:4.4}
0\le (k+1)Mt^{-1}f_k(|x|)\le\frac{1}{2},
\quad
z_*(x,t)\ge\zeta(t),
\end{equation}
for $x\in B(0,\delta\sqrt{t})$ and $t>T$.
Here the constant $\delta$ depends on $M$ and it is independent of $k$.
\end{lemma}
{\bf Proof.}
We see that $f_k$ is nonnegative and it satisfies 
$$
\frac{1}{\nu_k}\mbox{div}\,(\nu_k\nabla f_k)=1\quad\mbox{in}\quad {\bf R}^N.
$$
Furthermore, $\zeta'\le 0$ in $(T,\infty)$. 
Then we observe that $z_*(x,t)\le 2\zeta(t)$ and 
\begin{equation*}
\begin{split}
 & \partial_t z_*-\displaystyle{\frac{1}{\nu_k}\mbox{div}}\,(\nu_k\nabla z_*)\\
 & =2\zeta'(t)\left[1-(k+1)Mt^{-1}f_k(x)\right]
+2(k+1)M\zeta(t) t^{-2}f_k(x)
+2(k+1)Mt^{-1}\zeta(t)\\
& \ge 2\zeta'(t)+2(k+1)M t^{-1}\zeta(t)\ge 0
\end{split}
\end{equation*}
for $(x,t)\in{\bf R}^N\times(T,\infty)$. 
Furthermore, 
by Proposition~\ref{Proposition:2.4}~(b) 
we find $C>0$, independent of $k$, such that
$$
M(k+1)f_k(x)\le CM\int_0^{|x|}s\,ds=\frac{CM}{2}|x|^2
\le\frac{CM}{2}\delta^2t
$$
for $x\in B(0,\delta\sqrt{t})$ and $t>T$. 
Then, taking a sufficiently small $\delta\in(0,1]$, 
we see that $z_*$ satisfies \eqref{eq:4.4}. 
Thus Lemma~\ref{Lemma:4.2} follows.
$\Box$\vspace{3pt}

We are ready to prove Proposition~\ref{Proposition:4.1}.
\vspace{5pt}
\newline
{\bf Proof of Proposition~\ref{Proposition:4.1}.}  
The proof is divided into three steps. Let $j\in\{0,1,2,\dots\}$ and fix it.
%\vspace{3pt}
\newline
\underline{Step 1}: We prove inequality~\eqref{eq:4.1} for $0<t\le 1$. 
In this step the letter $C$ denotes generic positive constants independent of $k$. 
Let $\epsilon\in(0,1)$. 
Set 
$$
\tilde{v}_j(x,t):=\partial_t^j v_{k,i}(x,t),
\qquad
\tilde{w}_j(x,t):=\frac{\tilde{v}_j(x,t)}{h_k(|x|)},
\qquad
\eta_j(t):=t^{-\frac{N}{2p}-j-\frac{A_{1,k}}{2}},
$$
for $(x,t)\in{\bf R}^N\times(0,\infty)$. 
Then $\tilde{w}_j$ satisfies
\begin{equation}
\tag{$\mbox{W}_k$}
\partial_tw=\frac{1}{\nu_k}\mbox{div}\,(\nu_k\nabla w)\quad\mbox{in}\quad{\bf R}^N\times(0,\infty). 
\end{equation}
(See also Definition~\ref{Definition:2.1}.)
Since $\tilde{w}_j(\cdot,t)$ is radially symmetric for any $t>0$, 
 it follows from \cite[Lemma~2.2]{IM} that $\tilde{w}_j\in C({\bf R}^N\times(0,\infty))$. 
On the other hand, 
it follows from $h_0\in L^{p,\sigma}(B(0,1))$, \eqref{eq:1.3} and \eqref{eq:1.4} that $N+pA_{1,0}\ge 0$ 
and 
$$
A_{1,k}=\sqrt{\omega_k}(1+o(1))=k(1+o(1))\quad\mbox{as}\quad k\to\infty.
$$
Then we find $M_1>0$ such that 
$$
-M_1(k+1)\le -\frac{N}{2p}-j-\frac{A_{1,k}}{2}-\epsilon
\le -\frac{1}{2p}(N+pA_{1,0})-\epsilon<0
$$
for $k\in\{0,1,2,\dots\}$. 
This implies that 
$$
0\le -\frac{t}{t^{-\epsilon}\eta_j(t)}\frac{d}{dt}(t^{-\epsilon}\eta_j(t))\le M_1(k+1)\quad\mbox{for}\quad t>0. 
$$
Let $z_*$ and $\delta\in(0,1]$ be as in Lemma~\ref{Lemma:4.2} 
with $\zeta(t)$ and $T$ replaced by $t^{-\epsilon}\eta_j(t)$ and $0$, respectively. 
Here $\delta$ depends only on $M_1$ and it is independent of $k$ and $\epsilon$.

On the other hand, 
by Proposition~\ref{Proposition:2.3} 
we see that 
$$
h_k(|x|)\asymp |x|^{A_{1,k}}
\quad\mbox{for $x\in B(0,1)$ and $k\in\{0,1,2,\dots\}$}. 
$$
Then, by \eqref{eq:2.4} we have
\begin{equation}
\label{eq:4.5}
C^{-1}\eta_j(t)\le
t^{-\frac{N}{2}-j}\frac{\|h_0\|_{L^{p',\sigma'}(B(0,\sqrt{t}))}}{h_k(\delta\sqrt{t})h_0(\sqrt{t})}\le C\eta_j(t)
\quad\mbox{for}\quad t\in(0,1].
\end{equation}
This together with Lemmas~\ref{Lemma:4.1} and \ref{Lemma:4.2} implies that
\begin{equation}
\label{eq:4.6}
\begin{split}
|\tilde{w}_j(x,t+t_0)| & \le CM_{k,i}t^{-\frac{N}{2}-j}\frac{\|h_0\|_{L^{p',\sigma'}(B(0,\sqrt{t}))}}{h_k(\delta\sqrt{t})h_0(\sqrt{t})}\|\phi\|_{L^{p,\sigma}}\\
 & \le CM_{k,i}t^{-\epsilon}\eta_j(t)\|\phi\|_{L^{p,\sigma}}
 \le CM_{k,i}\|\phi\|_{L^{p,\sigma}}z_*(x,t)
\end{split}
\end{equation}
for $x\in\partial B(0,\delta\sqrt{t})$, $t\in(0,1]$, $t_0\in(0,1]$ and $k\in\{0,1,2,\dots\}$. 
Since $\tilde{w}_j\in C({\bf R}^N\times[0,\infty))$ and 
$$
\inf_{x\in B(0,\delta\sqrt{t})}z_*(x,t)\ge t^{-\epsilon}\eta_j(t)\to\infty\quad\mbox{as}\quad t\to +0,
$$
we see that 
\begin{equation}
\label{eq:4.7}
\tilde{w}_j(x,t_*+t_0)\le CM_{k,i}\|\phi\|_{L^{p,\sigma}}z_*(x,t_*)
\end{equation}
for $x\in B(0,\delta\sqrt{t})$ and small enough $t_*>0$. 
By Lemma~\ref{Lemma:4.2}, \eqref{eq:4.5}, \eqref{eq:4.6} and \eqref{eq:4.7}
we apply the comparison principle to obtain 
\begin{equation}
\label{eq:4.8}
\begin{split}
|\tilde{w}_j(x,t+t_0)| & \le CM_{k,i}\|\phi\|_{L^{p,\sigma}}z_*(x,t)
\le CM_{k,i}\|\phi\|_{L^{p,\sigma}}t^{-\epsilon}\eta_j(t)\\
 & \le CM_{k,i}
t^{-\frac{N}{2}-j-\epsilon}\frac{\|h_0\|_{L^{p',\sigma'}(B(0,\sqrt{t}))}}{h_k(\delta\sqrt{t})h_0(\sqrt{t})}\|\phi\|_{L^{p,\sigma}}
\end{split}
\end{equation}
for $x\in B(0,\delta\sqrt{t})$ and $t\in[t_*,1]$. 
Letting $t_*\to+0$, we see that \eqref{eq:4.8} holds for $x\in B(0,\delta\sqrt{t})$ and $t\in(0,1]$.
Since $\epsilon\in(0,1)$ and  $t_0\in(0,1)$ are arbitrary, 
we obtain inequality~\eqref{eq:4.1} for $x\in B(0,\delta\sqrt{t})$, $0<t\le 1$ and $k\in\{0,1,2,\dots\}$. 
It remains to prove inequality~\eqref{eq:4.1} for $t>1$.
\vspace{5pt}
\newline
\underline{Step 2}: 
Let $k\in\{0,1,2,\dots\}$ and fix it. 
We prove that inequality~\eqref{eq:4.1} holds for $t>1$. 
In this step the letter $C$ denotes generic positive constants depending on $k$ possibly. 
By \eqref{eq:4.8} we have
\begin{equation}
\label{eq:4.9}
|\tilde{w}_j(x,1)|\le CM_{k,i}\|\phi\|_{L^{p,\sigma}}\quad\mbox{for}\quad x\in B(0,\delta).
\end{equation}
Furthermore, by Lemma~\ref{Lemma:4.1} 
we see that 
\begin{equation}
\label{eq:4.10}
|\tilde{w}_j(x,t)|\le CM_{k,i}t^{-\frac{N}{2}-j}\frac{\|h_0\|_{L^{p',\sigma'}(B(0,\sqrt{t}))}}{h_k(\delta\sqrt{t})h_0(\sqrt{t})}\|\phi\|_{L^{p,\sigma}}
\end{equation}
for $x\in\partial B(0,\delta\sqrt{t})$ and $t\ge 1$. 
On the other hand, 
by \eqref{eq:1.5} we find $\alpha\in{\bf R}$ and $\beta\in{\bf R}$ such that 
\begin{equation}
\label{eq:4.11}
t^{-\frac{N}{2}-j}\frac{\|h_0\|_{L^{p',\sigma'}(B(0,\sqrt{t}))}}{h_k(\delta\sqrt{t})h_0(\sqrt{t})}
\asymp t^{-\alpha}[\log(1+t)]^\beta\quad\mbox{for}\quad t\ge 1.
\end{equation}

Assume either $\alpha> 0$ or $\alpha=0$, $\beta<0$.
Then we find $L\ge 1$ such that 
$$
\eta_2(t):=t^{-\alpha}[\log(L+t)]^\beta
$$
is monotone decreasing in $(1,\infty)$. 
Furthermore, we also find $M_2>0$ such that 
$$
0\le -\frac{t\eta'_2(t)}{\eta_2(t)}\le M_2\quad\mbox{in}\quad(1,\infty). 
$$
By Lemma~\ref{Lemma:4.2} 
we find a supersolution~$z_*$ to problem~$(\mbox{W}_k)$ and $\delta_k\in(0,\delta)$ such that 
\begin{equation}
\label{eq:4.12}
\begin{array}{ll}
 z_*\le 2\eta_2(t) & \mbox{for $(x,t)\in{\bf R}^N\times(1,\infty)$},\vspace{5pt}\\
 z_*(x,t)\ge\eta_2(t) & \mbox{for $(x,t)\in B(0,\delta_k\sqrt{t})\times(1,\infty)$}.
\end{array}
\end{equation}
Furthermore, by \eqref{eq:4.9} we have 
\begin{equation}
\label{eq:4.13}
|\tilde{w}_j(x,1)|
\le CM_{k,i}\|\phi\|_{L^{p,\sigma}}=\frac{CM_{k,i}}{\eta_2(1)}\|\phi\|_{L^{p,\sigma}}\eta_2(1)\le \frac{CM_{k,i}}{\eta_2(1)}\|\phi\|_{L^{p,\sigma}}z_*(x,1)
\end{equation}
for $x\in B(0,\delta_k)$. 
It follows from \eqref{eq:4.10}, \eqref{eq:4.11} and \eqref{eq:4.12} that
\begin{equation}
\label{eq:4.14}
|\tilde{w}_j(x,t)|\le CM_{k,i}\eta_2(t)\|\phi\|_{L^{p,\sigma}}
\le CM_{k,i}\|\phi\|_{L^{p,\sigma}}z_*(x,t)
\end{equation}
for $x\in\partial B(0,\delta_k\sqrt{t})$ and $t\ge 1$. 
By \eqref{eq:4.13} and \eqref{eq:4.14}, 
applying the comparison principle, 
we see that 
$$
|\tilde{w}_j(x,t)|\le CM_{k,i}\|\phi\|_{L^{p,\sigma}}z_*(x,t)
\quad\mbox{for $x\in B(0,\delta_k\sqrt{t})$ and $t\ge 1$}.
$$
This together with \eqref{eq:4.11} and \eqref{eq:4.12} implies that 
$$
|\tilde{w}_j(x,t)|
\le CM_{k,i}\|\phi\|_{L^{p,\sigma}}\eta_2(t)
\le CM_{k,i}t^{-\frac{N}{2}-j}\frac{\|h_0\|_{L^{p',\sigma'}(B(0,\sqrt{t}))}}{h_k(\delta\sqrt{t})h_0(\sqrt{t})}\|\phi\|_{L^{p,\sigma}}
$$
for $x\in B(0,\delta_k\sqrt{t})$ and $t\ge 1$.  
This implies that inequality~\eqref{eq:4.1} holds 
for $x\in B(0,\delta_k\sqrt{t})$ and $t\ge 1$
in the cases when $\alpha>0$ and when $\alpha=0$ and $\beta<0$. 
 
Assume either $\alpha<0$ or $\alpha=0$, $\beta\ge0$. 
Let $L'\ge 1$ be such that 
$$
\tilde{\eta}(t):=t^{-\alpha}[\log(L'+t)]^\beta
$$
is monotone increasing in $(1,\infty)$. 
Then $\tilde{\eta}$ is a supersolution to problem~$(\mbox{W}_k)$. 
Furthermore, similarly to \eqref{eq:4.13} and \eqref{eq:4.14}, 
by \eqref{eq:4.9}, \eqref{eq:4.10} and \eqref{eq:4.11} 
we have
\begin{equation*}
\begin{split}
|\tilde{w}_j(x,1)| & \le CM_{k,i}\|\phi\|_{L^{p,\sigma}}\tilde{\eta}(1)\quad\mbox{for}\quad x\in B(0,\delta_k),\\
|\tilde{w}_j(x,t)| & \le CM_{k,i}\|\phi\|_{L^{p,\sigma}}\tilde{\eta}(t)\quad\mbox{for}\quad
(x,t)\in \partial B(0,\delta_k\sqrt{t})\times(1,\infty).
\end{split}
\end{equation*}
Applying the comparison principle, we see that 
$|\tilde{w}_j(x,t)|\le CM_{k,i}\|\phi\|_{L^{p,\sigma}}\tilde{\eta}(t)$
for $(x,t)\in B(0,\delta_k\sqrt{t})\times[1,\infty)$. 
This together with \eqref{eq:4.11} implies that inequality~\eqref{eq:4.1} holds 
for $x\in B(0,\delta_k\sqrt{t})$ and $t\ge 1$
in the cases when $\alpha<0$ and when $\alpha=0$, $\beta<0$. 
Thus inequality~\eqref{eq:4.1} holds for $x\in B(0,\delta_k\sqrt{t})$ and $t\ge 1$ 
for fixed $k\in\{0,1,2,\dots\}$. 
\vspace{5pt}
\newline
\underline{Step 3}: 
We complete the proof of Proposition~\ref{Proposition:4.1}. 
Thanks to Steps 1 and 2, it suffices to prove inequality~\eqref{eq:4.1} 
for $x\in B(0,\delta\sqrt{t})$, $t\ge 1$ and large enough $k$. 
In this step the letter $C$ denotes generic positive constants independent of $k$. 

Let $k_*\in\{0,1,2,\dots\}$ be as in Proposition~\ref{Proposition:2.2}~(b) and $k\in\{k_*,k_*+1,\dots\}$.
Then $h_k$ is monotone increasing in $(1,\infty)$. 
We construct a supersolution to problem~$(\mbox{W}_k)$. 
Since 
$$
\frac{d}{dt}\left[h_k(\sqrt{t})^{-1}\right]
=-\frac{1}{2}t^{-\frac{1}{2}}h_k(\sqrt{t})^{-2}\left(\frac{d}{dr}h_k\right)(\sqrt{t}), 
$$
by Proposition~\ref{Proposition:2.2}~(b)
we find $M_3>0$ such that 
\begin{equation}
\label{eq:Z}
-kM_3^{-1} t^{-1}h_k(\sqrt{t})^{-1}
\ge \frac{d}{dt}h_k(\sqrt{t})^{-1}
\ge -kM_3 t^{-1}h_k(\sqrt{t})^{-1}
\end{equation}
for $t>0$ and $k\ge k_*$. 
Let $\tilde{\alpha}$, $\tilde{\beta}\in{\bf R}$ be such that 
\begin{equation}
\label{eq:4.11a}
t^{-\frac{N}{2}-j}\frac{\|h_0\|_{L^{p',\sigma'}(B(0,\sqrt{t}))}}{h_0(\sqrt{t})}
\asymp t^{-\tilde{\alpha}}[\log(1+t)]^{\tilde{\beta}}\quad\mbox{for}\quad t\ge 1.
\end{equation} 
Set 
$$
M_4:=|\tilde{\alpha}|+|\tilde{\beta}|(\log 2)^{-1}+M_3.
$$
By Lemma~\ref{Lemma:4.2}
we take small enough $\tilde{\delta}\in(0,\delta]$ to obtain 
$$
0\le (k+1)M_4t^{-1}f_k(|x|)\le\frac{1}{2}
$$
for $x\in B(0,\tilde{\delta}\sqrt{t})$ and $t>0$, where $f_k$ is as in Lemma~\ref{Lemma:4.2}. 
Set 
$$
\zeta_k(t):=t^{-\tilde{\alpha}}[\log(2+t)]^{\tilde{\beta}} h_k(\tilde{\delta}\sqrt{t})^{-1}\quad\mbox{for}\quad t>0. 
$$
Then, by \eqref{eq:Z} we have
\begin{equation}
\label{eq:4.15}
\begin{split}
 & \frac{d}{dt}\zeta_k(t)=-\tilde{\alpha} t^{-1}\zeta_k(t)+\tilde{\beta}[\log(2+t)]^{-1}(2+t)^{-1}\zeta_k(t)\\
 & \qquad\qquad\qquad
 +\tilde{\delta}^2 t^{-\tilde{\alpha}}[\log(2+t)]^{\tilde{\beta}}\frac{d}{d\tau}\left[h_k(\sqrt{\tau})^{-1}\right]\biggr|_{\tau=\tilde{\delta}^2 t}\\
 & \ge -\left[|\tilde{\alpha}|+|\tilde{\beta}|(\log 2)^{-1}\right]t^{-1}\zeta_k(t)-kM_3t^{-\tilde{\alpha}-1}[\log(2+t)]^{\tilde{\beta}} h_k(\tilde{\delta}\sqrt{t})^{-1}
 \ge -kM_4t^{-1}\zeta_k(t)
\end{split}
\end{equation} 
for $t>0$ and $k\ge k_*$. 
Similarly, taking sufficiently large $k_*$ if necessary, we see that 
\begin{equation}
\label{eq:4.16}
\frac{d}{dt}\zeta_k(t)\le 0\quad\mbox{for}\quad t>0. 
\end{equation}
Set 
$$
z_*(x,t):=2\zeta_k(t)\left[1-(k+1)M_4t^{-1}f_k(|x|)\right].
$$
By \eqref{eq:4.15} and \eqref{eq:4.16}, 
taking small enough $\tilde{\delta}$ if necessary, 
we apply Lemma~\ref{Lemma:4.2} to see that $z_*$ is a supersolution to problem~$(\mbox{W}_k)$
and 
\begin{equation}
\label{eq:4.17}
2\zeta_k(t)\ge z_*(x,t)\ge\zeta_k(t)\quad\mbox{for}\quad x\in B(0,\tilde{\delta}\sqrt{t})\quad\mbox{and}\quad t>0. 
\end{equation}
Here $\tilde{\delta}\in(0,\delta]$ depends only on $M_1$ and $M_4$ and it is independent of $k$. 
On the other hand, 
by~\eqref{eq:4.8} with $\delta=\tilde{\delta}$ we have
\begin{equation}
\label{eq:4.18}
|\tilde{w}_j(x,1)|\le CM_{k,i}\|\phi\|_{L^{p,\sigma}}\frac{\zeta_k(1)}{h_k(\tilde{\delta})\zeta_k(1)}
\le CM_{k,i}\|\phi\|_{L^{p,\sigma}}\zeta_k(1)
\le CM_{k,i}\|\phi\|_{L^{p,\sigma}}z_*(x,1)
\end{equation}
for $x\in B(0,\tilde{\delta})$.
On the other hand, by Lemma~\ref{Lemma:4.1} and \eqref{eq:4.11a}
we obtain 
\begin{equation}
\label{eq:4.19}
\begin{split}
|\tilde{w}_j(x,t)|
\le CM_{k,i}t^{-\frac{N}{2}-j}\frac{\|h_0\|_{L^{p',\sigma'}(B(0,\sqrt{t}))}}{h_k(\tilde{\delta}\sqrt{t})h_0(\sqrt{t})}\|\phi\|_{L^{p,\sigma}}
\le CM_{k,i}\|\phi\|_{L^{p,\sigma}}\zeta(t)
\end{split}
\end{equation}
for $x\in\partial B(0,\tilde{\delta}\sqrt{t})$ and $t>0$. 
By \eqref{eq:4.17}, \eqref{eq:4.18} and \eqref{eq:4.19} 
we apply the comparison principle to obtain 
$$
|\tilde{w}_j(x,t)|\le CM_{k,i}\|\phi\|_{L^{p,\sigma}}z_*(x,t)
\le CM_{k,i}\|\phi\|_{L^{p,\sigma}}\zeta(t)
$$
for $x\in B(0,\tilde{\delta}\sqrt{t})$ and $t>0$. 
This together with \eqref{eq:4.11a} implies that 
inequality~\eqref{eq:4.1} holds 
for $x\in B(0,\tilde{\delta}\sqrt{t})$, $t\ge 1$ and $k\ge k_*$.
Thus Proposition~\ref{Proposition:4.1} follows.
$\Box$
%%%%%%%%%%%%%%%%%%%%%%%%%%%
\subsection{Estimates of spacial derivatives}
%%%%%%%%%%%%%%%%%%%%%%%%%%%
In this subsection, thanks to Proposition~\ref{Proposition:4.1}, 
we obtain the following proposition.
\begin{proposition}
\label{Proposition:4.2}
Assume the same conditions as in Proposition~{\rm\ref{Proposition:4.1}}. 
Set 
$$
w_{k,i}(|x|,t):=\frac{v_{k,i}(|x|,t)}{h_k(|x|)}\quad\mbox{for}\quad(x,t)\in{\bf R}^N\times(0,\infty),
$$
where $k\in\{0,1,2,\dots\}$ and $i\in\{1,\dots,d_k\}$.
Then, for any $j\in\{0,1,2,\dots\}$,
\begin{equation}
\label{eq:4.20}
\partial_t^j w_{k,i}(|x|,t)=\partial_t^j w_{k,i}(0,t)+F_{k,i}^j(x,t),
\quad (x,t)\in{\bf R}^N\times(0,\infty),
\end{equation}
where 
$$
F_{k,i}^j(x,t):=\int_0^{|x|}s^{1-N}\nu_k(s)^{-1}\left(\int_0^s \tau^{N-1}\nu_k(\tau)[\partial_t^{j+1}w_{k,i}](\tau,t)\,d\tau\right)\,ds.
$$
Furthermore, there exists $C>0$ and $\delta\in(0,1]$ such that 
\begin{eqnarray}
\label{eq:A}
 & & |\partial_t^j w_{k,i}(x,t)|\le CM_{k,i}t^{-\frac{N}{2}-j}\frac{\|h_0\|_{L^{p',\sigma'}(B(0,\sqrt{t}))}}{h_k(\delta\sqrt{t})h_0(\sqrt{t})}\|\phi\|_{L^{p,\sigma}},
\qquad\qquad\\
\label{eq:B}
 & & 
 \left|\nabla^\ell F_{k,i}^j(x,t)\right|
 \le CM_{k,i}t^{-\frac{N}{2}-j-1}|x|^{2-\ell}\frac{\|h_0\|_{L^{p',\sigma'}(B(0,\sqrt{t}))}}{h_k(\delta\sqrt{t})h_0(\sqrt{t})}\|\phi\|_{L^{p,\sigma}},
\end{eqnarray}
for $x\in B(0,\delta\sqrt{t})$ and $t>0$, where $\ell\in\{0,1\}$. 
\end{proposition}
{\bf Proof.}
Let $j\in\{0,1,2,\dots\}$. 
By Proposition~\ref{Proposition:4.1} 
we find $C>0$ and $\delta\in(0,1]$ such that 
\begin{equation}
\label{eq:C}
|\partial_t^{j+1} w_{k,i}(x,t)|\le CM_{k,i}t^{-\frac{N}{2}-j-1}\frac{\|h_0\|_{L^{p',\sigma'}(B(0,\sqrt{t}))}}{h_k(\delta\sqrt{t})h_0(\sqrt{t})}\|\phi\|_{L^{p,\sigma}}
\end{equation}
for $x\in B(0,\delta\sqrt{t})$, $t>0$ and $k\in\{0,1,2,\dots\}$.
This implies \eqref{eq:A}.
Furthermore, combining \eqref{eq:C} and Proposition~\ref{Proposition:2.4}~(b), we have
\begin{equation}
\label{eq:4.21}
\begin{split}
 & \left||x|^{1-N}\nu_k(|x|)^{-1}\int_0^{|x|} \tau^{N-1}\nu_k(\tau)(\partial_t^{j+1} w_{k,i})(\tau,t)\,d\tau\right|\\
 & \qquad
 \le \|\partial_t^{j+1}w_{k,i}(t)\|_{L^\infty(B(0,\delta\sqrt{t}))}
 \cdot |x|^{1-N}\nu_k(|x|)^{-1}\int_0^{|x|} \tau^{N-1}\nu_k(\tau)\,d\tau\\
 & \qquad
 \le C(k+1)^{-1}|x|\|\partial_t^{j+1}w_{k,i}(t)\|_{L^\infty(B(0,\delta\sqrt{t}))},\\
 & |\nabla^\ell F_{k,i}^j(x,t)|
 \le C(k+1)^{-1}M_{k,i}t^{-\frac{N}{2}-j-1}|x|^{2-\ell}\frac{\|h_0\|_{L^{p',\sigma'}(B(0,\sqrt{t}))}}{h_k(\delta\sqrt{t})h_0(\sqrt{t})}\|\phi\|_{L^{p,\sigma}},
\end{split}
\end{equation}
for $x\in B(0,\delta\sqrt{t})$, $t>0$ and $k\in\{0,1,2,\dots\}$, where $\ell\in\{0,1\}$. 
Then $F_{k,i}^j$ is well-defined and \eqref{eq:B} holds. 
Furthermore, it satisfies 
$$
\frac{1}{\nu_k}\mbox{div}\,(\nu_k\nabla F_{k,i}^j)=\partial_t^{j+1}w_{k,i}
\quad\mbox{in}\quad{\bf R}^N\times(0,\infty). 
$$
Set 
$$
\tilde{w}_{k,i}^j(|x|,t):=\partial_t^j w_{k,i}(|x|,t)-F_{k,i}^j(|x|,t).
$$
Then it follows that 
\begin{equation}
\label{eq:4.22}
\frac{1}{\nu_k}\mbox{div}\,(\nu_k\nabla \tilde{w}_k)=0\quad\mbox{in}\quad{\bf R}^N\times(0,\infty).
\end{equation}
For any fixed $t>0$, 
set $z_{k,i}^j(r):=h_k(r)\tilde{w}_{k,i}^j(r,t)$ for $r>0$. 
Then, by \eqref{eq:4.22} we see that $z_{k,i}^j$ satisfies 
\begin{equation*}
\begin{split}
 & \frac{d^2}{dr^2}z_{k,i}^j+\frac{N-1}{r}\frac{d}{dr}z_{k,i}^j-V_k(r)z_{k,i}^j=0
\quad\mbox{for}\quad r\in(0,\infty),\\
 & \lim_{r\to +0}\frac{z_{k,i}^j(r)}{h_k(r)}=\tilde{w}_{k,i}^j(0,t)=\partial_t^j w_{k,i}(0,t).
\end{split}
\end{equation*}
By \eqref{eq:1.4} we see that $z_{k,i}^j(r)=\partial_t^j w_{k,i}(0,t)h_k(r)$ for $r\in(0,\infty)$, that is, 
$$
\partial_t^j w_{k,i}(|x|,t)-F_{k,i}^j(|x|,t)
=\tilde{w}_{k,i}^j(|x|,t)=\partial_t^j w_{k,i}(0,t),\quad x\in{\bf R}^N,
$$
for $t>0$. Thus relation~\eqref{eq:4.20} holds. 
Thus Proposition~\ref{Proposition:4.2} follows.
$\Box$
%%%%%%%%%%%%%%%%%%%%%%%%%%%%%%%%%%%%
%%%%%%%%%%%%%%%%%%%%%%%%%%%%%%%%%%%%
\section{Proof of Theorem~\ref{Theorem:1.1}}
%%%%%%%%%%%%%%%%%%%%%%%%%%%%%%%%%%%%
%%%%%%%%%%%%%%%%%%%%%%%%%%%%%%%%%%%%
We complete the proof of Theorem~\ref{Theorem:1.1}, and prove Corollary~\ref{Corollary:1.1}.
\vspace{5pt}
\newline
{\bf Proof of Theorem~\ref{Theorem:1.1}.}
Theorem~\ref{Theorem:1.1} with $\ell=0$ follows from Proposition~\ref{Proposition:3.1}~(b). 
It suffices to prove \eqref{eq:1.7} with $\ell=1$. 

Let $(p,q,\sigma,\theta)\in\Lambda$ and $\ell=1$. 
We can assume, without loss of generality, that $h_0\in L^{p',\sigma'}(B(0,1))$ and $\nabla h_0\in L^{q,\theta}(B(0,1))$. 
Then, by the definition of the Lorentz norm and Proposition~\ref{Proposition:2.1}
we see that $h_0\in L^{q,\theta}(B(0,1))$. These imply that
\begin{equation}
\label{eq:s}
h_0\in L^{p',\sigma'}(B(0,1))\cap L^{p,\sigma}(B(0,1))\cap L^{q,\theta}(B(0,1)). 
\end{equation}
Let $\phi\in C_c({\bf R}^N)$. 
We use the same notations as in \eqref{eq:1.8}, \eqref{eq:1.9} and \eqref{eq:1.10}.
By \eqref{eq:1.1} we apply the regularity theorems for elliptic equations to obtain
$$
\|Q_{k,i}\|_{L^\infty({\bf S}^{N-1})}+
\|\nabla Q_{k,i}\|_{L^\infty({\bf S}^{N-1})}\le C(\omega_k+1)\|Q_{k,i}\|_{L^2({\bf S}^{N-1})}\le C(k+1)^2
$$
for $\ell\in\{0,1\}$, $k\in\{0,1,2,\dots\}$ and $i\in\{1,\dots,d_k\}$. 
Since 
$$
u_{k,i}(x,t)=w_{k,i}(x,t)h_k(|x|)Q_{k,i}\left(\frac{x}{|x|}\right),
$$
applying Proposition~\ref{Proposition:4.2} with \eqref{eq:s}, we obtain
\begin{equation}
\label{eq:5.1}
\begin{split}
|\partial_t^j\nabla u_{k,i}(x,t)| 
 & \le C(k+1)^4t^{-\frac{N}{2}-j}\frac{\|h_0\|_{L^{p',\sigma'}(B(0,\sqrt{t}))}}{h_k(\delta\sqrt{t})h_0(\sqrt{t})}\|\phi\|_{L^{p,\sigma}}\\
 & \qquad\times\left\{t^{-1}|x|h_k(x)
+|\nabla h_k(x)|+|x|^{-1}h_k(x)\right\}
\end{split}
\end{equation}
for $x\in B(0,\delta\sqrt{t})\setminus\{0\}$ and $t>0$.
Here $\delta$ is as in Proposition~\ref{Proposition:4.2}.
By Propositions~\ref{Proposition:2.1}~and~\ref{Proposition:2.2} 
we see that
\begin{equation}
\label{eq:5.2}
   |\partial_t^j\nabla u_{k,i}(x,t)|\le C(k+1)^5
   t^{-\frac{N}{2}-j}\frac{\|h_0\|_{L^{p',\sigma'}(B(0,\sqrt{t}))}}{h_k(\delta\sqrt{t})h_0(\sqrt{t})}\|\phi\|_{L^{p,\sigma}}|x|^{-1}h_k(x)
\end{equation}
for $x\in B(0,\delta\sqrt{t})\setminus\{0\}$ and $t>0$.
Set 
$$
\tilde{u}(x,t):=u(x,t)-u_{0,1}(x,t)=\sum_{k=1}^\infty\sum_{i=1}^{d_k}u_{k,i}(x,t).
$$
By \eqref{eq:1.2}, \eqref{eq:1.10} and \eqref{eq:5.2} we see that
\begin{equation}
\label{eq:5.3}
\begin{split}
 & |\partial_t^j\nabla \tilde{u}(x,t)|\le\sum_{k=1}^\infty\sum_{i=1}^{d_k}|\partial_t^j\nabla u_{k,i}(x,t)|\\
 & \qquad
\le C\|\phi\|_{L^{p,\sigma}}t^{-\frac{N}{2}-j}|x|^{-1}\frac{\|h_0\|_{L^{p',\sigma'}(B(0,\sqrt{t}))}}{h_0(\sqrt{t})}
\sum_{k=1}^\infty\sum_{i=1}^{d_k}(k+1)^5\frac{h_k(|x|)}{h_k(\delta\sqrt{t})}\\
 & \qquad
 \le C\|\phi\|_{L^{p,\sigma}}t^{-\frac{N}{2}-j}|x|^{-1}\frac{\|h_0\|_{L^{p',\sigma'}(B(0,\sqrt{t}))}}{h_0(\sqrt{t})}
\sum_{k=1}^\infty (k+1)^{N+3}\frac{h_k(|x|)}{h_k(\delta\sqrt{t})}
\end{split}
\end{equation}
for $x\in B(0,\delta\sqrt{t})\setminus\{0\}$ and $t>0$.

For $k\in\{1,2,\dots\}$, set 
$$
\iota_k(r):=\left\{
\begin{array}{ll}
r^{A_{1,k}-A_{1,1}} & \mbox{for}\quad \mbox{$0<r<1$},\vspace{5pt}\\
r^{A_{2,k}-A_{2,1}} & \mbox{for}\quad \mbox{$r\ge 1$}.
\end{array}
\right.
$$
Then $\iota_k$ is monotone increasing in $(0,\infty)$. 
Furthermore, by Proposition~\ref{Proposition:2.3} we have 
$$
\frac{h_k(r)}{h_1(r)}\asymp\iota_k(r)\quad\mbox{for $r>0$ and $k\in\{1,2,\dots\}$}. 
$$
Since $A_{i,k}-A_{i,1}=k(1+o(1))$ as $k\to\infty$, where $i=1,2$, 
we find $\gamma>0$ such that 
$$
A_{i,k}-A_{i,1}\ge\frac{k}{2}-\gamma\quad\mbox{for}\quad k\in\{1,2,\dots\}.
$$
Then we see that $\iota_k(\epsilon r)\le\epsilon^{\frac{k}{2}-\gamma}\iota_k(r)$ for $r>0$ and $\epsilon\in(0,1)$. 
This implies that
\begin{equation*}
\begin{split}
\frac{h_k(|x|)}{h_k(\delta\sqrt{t})}
 & =\frac{h_1(|x|)}{h_1(\delta\sqrt{t})}\frac{h_k(|x|)}{h_1(|x|)}\frac{h_1(\delta\sqrt{t})}{h_k(\delta\sqrt{t})}\\
 & \le C\frac{h_1(|x|)}{h_1(\delta\sqrt{t})}\frac{\iota_k(|x|)}{\iota_k(\delta\sqrt{t})}
 \le C\frac{h_1(|x|)}{h_1(\delta\sqrt{t})}\frac{\iota_k(\epsilon\delta\sqrt{t})}{\iota_k(\delta\sqrt{t})}
 \le C\epsilon^{\frac{k}{2}-\gamma}\frac{h_1(|x|)}{h_1(\delta\sqrt{t})}
\end{split}
\end{equation*}
for $x\in B(0,\epsilon\delta\sqrt{t})\setminus\{0\}$, $t>0$, $\epsilon\in(0,1)$ and $k\in\{1,2,\dots\}$. 
Taking small enough $\epsilon>0$ if necessary, we see that
\begin{equation}
\label{eq:5.4}
\sum_{k=1}^\infty(k+1)^{N+3}\frac{h_k(|x|)}{h_k(\delta\sqrt{t})}
\le C\frac{h_1(|x|)}{h_1(\delta\sqrt{t})}\sum_{k=1}^\infty \epsilon^{\frac{k}{2}-\gamma}(k+1)^{N+3}
\le C_\epsilon\frac{h_1(|x|)}{h_1(\delta\sqrt{t})}
\end{equation}
for $x\in B(0,\epsilon\delta\sqrt{t}))\setminus\{0\}$ and $t>0$. Here $C_\epsilon$ is a positive constant depending on $\epsilon>0$.
By \eqref{eq:5.3} and \eqref{eq:5.4} we obtain 
\begin{equation}
\label{eq:5.5}
|\partial_t^j\nabla \tilde{u}(x,t)|
\le CC_\epsilon\|\phi\|_{L^{p,\sigma}}t^{-\frac{N}{2}-j}|x|^{-1}\frac{\|h_0\|_{L^{p',\sigma'}(B(0,\sqrt{t}))}}{h_0(\sqrt{t})}\frac{h_1(|x|)}{h_1(\delta\sqrt{t})}
\end{equation}
for $x\in B(0,\epsilon\delta\sqrt{t})\setminus\{0\}$ and $t>0$. 
On the other hand, 
by \eqref{eq:5.1} and \eqref{eq:5.2} we have 
\begin{equation}
\label{eq:5.6}
\begin{split}
|\partial_t^j\nabla u_{0,1}(x,t)| & 
\le Ct^{-\frac{N}{2}-j-1}|x|\frac{\|h_0\|_{L^{p',\sigma'}(B(0,\sqrt{t}))}}{h_0(\sqrt{t})h_0(\delta\sqrt{t})}\|\phi\|_{L^{p,\sigma}}|h_0(x)|\\
 & +Ct^{-\frac{N}{2}-j}\frac{\|h_0\|_{L^{p',\sigma'}(B(0,\sqrt{t}))}}{h_0(\sqrt{t})h_0(\delta\sqrt{t})}\|\phi\|_{L^{p,\sigma}}|\nabla h_0(x)|\\
 & \le Ct^{-\frac{N}{2}-j}\frac{\|h_0\|_{L^{p',\sigma'}(B(0,\sqrt{t}))}}{h_0(\sqrt{t})h_0(\delta\sqrt{t})}\|\phi\|_{L^{p,\sigma}}
 \left[|\nabla h_0(x)|+t^{-\frac{1}{2}}h_0(|x|)\right]
\end{split}
\end{equation}
for $x\in B(0,\epsilon\delta\sqrt{t})\setminus\{0\}$ and $t>0$. 
Therefore, 
combining Proposition~\ref{Proposition:3.1}~(a),  \eqref{eq:5.5} and \eqref{eq:5.6}, 
we obtain 
\begin{equation*}
\begin{split}
 & \|\partial_t^j\nabla u(t)\|_{L^{q,\theta}}\\
 & \le \|\partial_t^j\nabla u_{0,1}(t)\|_{L^{q,\theta}(B(0,\epsilon\delta\sqrt{t}))}+\|\partial_t^j\nabla\tilde{u}(t)\|_{L^{q,\theta}(B(0,\epsilon\delta\sqrt{t}))}
+\|\partial_t^j\nabla u(t)\|_{L^{q,\theta}(B(0,\epsilon\delta\sqrt{t}))^c}\\
 & \le Ct^{-\frac{N}{2}-j}\frac{\|h_0\|_{L^{p',\sigma'}(B(0,\sqrt{t}))}}{h_0(\sqrt{t})}\|\phi\|_{L^{p,\sigma}}\\
 & \quad\times\left[\frac{\|\nabla h_0\|_{L^{q,\theta}(B(0,\epsilon\delta\sqrt{t}))}}{h_0(\delta\sqrt{t})}
 +t^{-\frac{1}{2}}\frac{\|h_0\|_{L^{q,\theta}(B(0,\epsilon\delta\sqrt{t}))}}{h_0(\delta\sqrt{t})}
 +\frac{\|\tilde{h}_1\|_{L^{q,\theta}(B(0,\epsilon\delta\sqrt{t}))}}{h_1(\delta\sqrt{t})}\right]
\end{split}
\end{equation*}
for $t>0$, where $\tilde{h}_1(x):=|x|^{-1}h_1(|x|)$. 
This together with \eqref{eq:1.4} and \eqref{eq:1.5} implies that 
\begin{equation}
\label{eq:5.7}
\begin{split}
\|\partial_t^j\nabla e^{-tH}\|_{(L^{p,\sigma}\to L^{q,\theta})}
 & \le Ct^{-\frac{N}{2}-j}\frac{\|h_0\|_{L^{p',\sigma'}(B(0,\sqrt{t}))}}{h_0(\sqrt{t})}\\
 & \hspace{-5pt}\times\left[\frac{\|\nabla h_0\|_{L^{q,\theta}(B(0,\sqrt{t}))}}{h_0(\sqrt{t})}
 +t^{-\frac{1}{2}}\frac{\|h_0\|_{L^{q,\theta}(B(0,\sqrt{t}))}}{h_0(\sqrt{t})}
 +\frac{\|\tilde{h}_1\|_{L^{q,\theta}(B(0,\sqrt{t}))}}{h_1(\sqrt{t})}\right]
\end{split}
\end{equation}
for $t>0$. 

On the other hand, 
by Proposition~\ref{Proposition:2.1}
we find $R>0$ such that 
\begin{equation}
\label{eq:5.8}
\begin{array}{ll}
h_0(|x|)\le C|x||\nabla h_0(|x|)|\le Ct^{\frac{1}{2}}|\nabla h_0(|x|)| & \mbox{if}\quad \lambda_1\not=0,\vspace{5pt}\\
h_0(|x|)\le C & \mbox{if}\quad \lambda_1=0,
\end{array}
\end{equation}
for $x\in [B(0,R)\cap B(0,\sqrt{t})]\setminus\{0\}$. 
Furthermore, taking small enough $R>0$ if necessary, 
by \eqref{eq:1.4} we see that 
$$
\frac{h_1(r)}{h_0(r)}\le C\frac{h_1(\sqrt{t})}{h_0(\sqrt{t})}
\quad\mbox{for}\quad x\in B(0,\sqrt{t})\cap B(0,R).
$$
These imply that 
\begin{equation}
\label{eq:5.9}
\begin{split}
\frac{\tilde{h}_1(|x|)}{h_1(\sqrt{t})} & =
\frac{|x|^{-1}h_1(|x|)}{h_1(\sqrt{t})}
\le C|x|^{-1}\frac{h_0(|x|)}{h_0(\sqrt{t})}
\le C\frac{|\nabla h_0(|x|)|}{h_0(\sqrt{t})}\quad\mbox{if}\quad\lambda_1\not=0,\\
\frac{\tilde{h}_1(|x|)}{h_1(\sqrt{t})} & \le C|x|^{-1}\frac{|x|}{\sqrt{t}}\le Ct^{-\frac{1}{2}}\quad\mbox{if}\quad\lambda_1=0,
\end{split}
\end{equation}
for $x\in[B(0,\sqrt{t})\cap B(0,R)]\setminus\{0\}$. 
By \eqref{eq:5.7}, \eqref{eq:5.8} and \eqref{eq:5.9} we obtain 
$$
\|\partial_t^j\nabla e^{-tH}\|_{(L^{p,\sigma}\to L^{q,\theta})}
 \le Ct^{-\frac{N}{2}-j}\frac{\|h_0\|_{L^{p',\sigma'}(B(0,\sqrt{t}))}}{h_0(\sqrt{t})}
\left[\frac{\|\nabla h_0\|_{L^{q,\theta}(B(0,\sqrt{t}))}}{h_0(\sqrt{t})}+t^{\frac{N}{2q}-\frac{1}{2}}\right]
$$
for $t\in(0,\sqrt{R})$, which implies \eqref{eq:1.7} for $t\in(0,\sqrt{R})$. 
Therefore Theorem~\ref{Theorem:1.1} follows in the case when $0<t<\sqrt{R}$. 

It remains to prove that \eqref{eq:1.7} holds for $t\ge \sqrt{R}$.
Similarly to \eqref{eq:5.8}, 
by Proposition~\ref{Proposition:2.2}
we find $R'\in(R,\infty)$ such that 
\begin{equation}
\label{eq:5.10}
\begin{array}{ll}
h_0(|x|)\le C & \mbox{if}\quad v_0\equiv1, \vspace{5pt}\\
h_0(|x|)\le C|x||\nabla h_0(|x|)|\le Ct^{\frac{1}{2}}|\nabla h_0(|x|)| & \mbox{otherwise},
\end{array}
\end{equation}
for $x\in B(0,R')^c\cap B(0,\sqrt{t})$. 
Furthermore, similarly to \eqref{eq:5.9},
taking large enough $R'$ if necessary, 
by \eqref{eq:1.5} 
we see that $h_1(r)/h_0(r)\asymp v_1(r)/v_0(r)$
and $v_1(r)/v_0(r)$ is monotone increasing in $(R',\infty)$ and obtain
\begin{equation}
\label{eq:5.11}
\begin{split}
\frac{\tilde{h}_1(|x|)}{h_1(\sqrt{t})} & \le C|x|^{-1}\frac{|x|}{\sqrt{t}}\le Ct^{-\frac{1}{2}}\quad\mbox{if}\quad v_0\equiv 1,\\
\frac{\tilde{h}_1(|x|)}{h_1(\sqrt{t})} & =
\frac{|x|^{-1}h_1(|x|)}{h_1(\sqrt{t})}
\le C|x|^{-1}\frac{h_0(|x|)}{h_0(\sqrt{t})}
\le C\frac{|\nabla h_0(|x|)|}{h_0(\sqrt{t})}\quad\mbox{otherwise},
\end{split}
\end{equation}
for $x\in B(0,R')^c\cap B(0,\sqrt{t})$. 
Therefore, by \eqref{eq:5.7}, \eqref{eq:5.10} and \eqref{eq:5.11} we obtain 
\begin{equation*}
 \|\partial_t^j\nabla e^{-tH}\|_{(L^{p,\sigma}\to L^{q,\theta})}
\le Ct^{-\frac{N}{2}-j}\frac{\|h_0\|_{L^{p',\sigma'}(B(0,\sqrt{t}))}}{h_0(\sqrt{t})}
 \left[\frac{\|\nabla h_0\|_{L^{q,\theta}(B(0,\sqrt{t}))}}{h_0(\sqrt{t})}+t^{\frac{N}{2q}-\frac{1}{2}}\right]
\end{equation*}
for $t\in[\sqrt{R},\infty)$. This implies \eqref{eq:1.7} for $t\in[\sqrt{R},\infty)$. 
Thus Theorem~\ref{Theorem:1.1} follows. 
$\Box$
\vspace{5pt}

\noindent
{\bf Proof of Corollary~\ref{Corollary:1.1}.}
Corollary~\ref{Corollary:1.1} easily follows from Theorem~\ref{Theorem:1.1} and Proposition~\ref{Proposition:2.4}~(a).
$\Box$
%%%%%%%%%%%%%%%%%%%%%%%%%%%%%%%%%%%%
%%%%%%%%%%%%%%%%%%%%%%%%%%%%%%%%%%%%
%\noindent
%{\bf Acknowledgements.} 
%The first author was supported in part by the Grant-in-Aid for Scientific Research (A)(No.~15H02058)
%from Japan Society for the Promotion of Science. 
%%%%%%%%%%%%%%%%%%%%%%%%%%%%%%%%%%%%%%
%%%%%%%%%%%%    references    %%%%%%%%%%%%%%%%%%
%%%%%%%%%%%%%%%%%%%%%%%%%%%%%%%%%%%%%%
  
%%%%%%%%%%%%%%%%%%%%%%%%%%%%%%%%%%%%
%%%%%%%%%%%%%%%%%%%%%%%%%%%%%%%%%%%% 
\end{document}